\documentclass[11pt,twoside]{article}

\usepackage{amsbsy,amsfonts,amsmath,amssymb,eucal,mathrsfs}
\usepackage[all]{xy}

\def\itemn#1{\item[\hspace{0.6mm} {\rm (#1)}]}
\renewcommand{\ge}{\geqslant}
\renewcommand{\le}{\leqslant}
\def\zmod#1{\mathbb{Z}/#1\mathbb{Z}}
\def\longto{\longrightarrow}

\renewcommand{\tilde}{\widetilde}

\renewcommand{\bar}{\overline}
\newcommand{\sLie}{\mathscr{L}ie}

\newcommand{\N}{\mathbb{N}}
\newcommand{\Z}{\mathbb Z}

\newcommand{\Q}{\mathbb{Q}}

\newcommand{\F}{\mathbb{F}}

\newcommand{\g}{\mathfrak g}
\newcommand{\h}{\mathfrak h}
\renewcommand{\u}{\mathfrak u}

\newcommand{\divise}{\ | \ }

\renewcommand{\iff}{if and only if }
\newcommand{\tr}{{}^t}

\newcommand{\scal}[1]{\langle #1 \rangle}

\newcommand{\lpara}{
\ \vspace{3pt}

\noindent}

\newcommand{\matdd}[4]{
\left (
\begin{array}{cc}
#1 & #2  \\
#3 & #4
\end{array}
\right )   }

\newcommand{\matddr}[4]{
\left (
\begin{array}{cc}
\hspace{-.2cm} #1 & \hspace{-.2cm}#2  \\
\hspace{-.2cm} #3 & \hspace{-.2cm}#4
\end{array}
\hspace{-.2cm} \right )   }

\newtheorem{definition}{Definition}[subsection]
\newenvironment{defi}{\begin{definition} \rm}{\end{definition}}

\newtheorem{prop}[definition]{Proposition}
\newtheorem{lemm}[definition]{Lemma}
\newtheorem{fait}[definition]{Fact}
\newtheorem{coro}[definition]{Corollary}
\newtheorem{theo}[definition]{Theorem}
\newtheorem{notation}[definition]{Notation}
\newenvironment{nota}{\begin{notation} \rm}{\end{notation}}
\newtheorem{remark}[definition]{Remark}
\newenvironment{rema}{\begin{remark} \rm}{\end{remark}}
\newtheorem{example}[definition]{Example}

\newtheorem{nothing}[definition]{$\!\!$}
\newenvironment{noth}{\begin{nothing} \rm}{\end{nothing}}
\newenvironment{proo}{{\flushleft \bf Proof :}}{\hfill $\square$ \vspace{5mm}}

\DeclareMathOperator{\Spec}{Spec}

\DeclareMathOperator{\End}{End}

\DeclareMathOperator{\Hom}{Hom}
\DeclareMathOperator{\Ext}{Ext}

\DeclareMathOperator{\Id}{Id}
\DeclareMathOperator{\pf}{pf}
\DeclareMathOperator{\Lie}{Lie}

\DeclareMathOperator{\pr}{pr}
\DeclareMathOperator{\ad}{ad}
\DeclareMathOperator{\Ad}{Ad}
\DeclareMathOperator{\Int}{Int}
\DeclareMathOperator{\res}{res}
\DeclareMathOperator{\coker}{coker}
\DeclareMathOperator{\Reg}{Reg}
\DeclareMathOperator{\Sing}{Sing}
\DeclareMathOperator{\Sym}{Sym}

   \def\cE{{\cal E}}
\def\cF{{\cal F}} \def\cG{{\cal G}} \def\cH{{\cal H}} \def\cI{{\cal I}}
\def\cK{{\cal K}}  \def\cM{{\cal M}} \def\cN{{\cal N}}
\def\cO{{\cal O}}

\def\bA{{\mathbb A}}  
\def\bF{{\mathbb F}}  
  \def\bP{{\mathbb P}}
\def\bQ{{\mathbb Q}}  \def\bV{{\mathbb V}}
 \def\bZ{{\mathbb Z}}

 \def\sM{{\sf M}}

\def\fb{\mathfrak{b}} \def\fg{\mathfrak{g}} \def\fh{\mathfrak{h}}
\def\fo{\mathfrak{o}} \def\ft{\mathfrak{t}} \def\fu{\mathfrak{u}}
\def\fz{\mathfrak{z}}

\def\fS{\mathfrak{S}}

\def\fgl{\mathfrak{gl}} \def\fsl{\mathfrak{sl}}
\def\fso{\mathfrak{so}} \def\fpso{\mathfrak{pso}}
\def\fsp{\mathfrak{sp}} \def\fpsp{\mathfrak{psp}}
\def\fpsl{\mathfrak{psl}} \def\fspin{\mathfrak{spin}}
\def\fpspin{\mathfrak{pspin}}

\begin{document}

\begin{center}
{\bf \Large On the adjoint quotient of Chevalley}

\medskip

{\bf \Large groups over arbitrary base schemes}

\bigskip
\bigskip

{\bf Pierre-Emmanuel Chaput and Matthieu Romagny}

\bigskip

\end{center}

\bigskip

{\def\thefootnote{\relax}
\footnote{ \hspace{-6.8mm}
}
}

\begin{center}
\bf{Abstract}
\end{center}

For a split semisimple Chevalley group scheme $G$ with Lie algebra $\g$ over
an arbitrary base scheme $S$, we consider the quotient of $\fg$ by the adjoint
action of $G$. We study in detail the structure of $\fg$ over $S$. Given a
maximal torus $T$ with Lie algebra $\ft$ and associated Weyl group~$W$, we
show that the Chevalley morphism $\pi:\ft/W\to \fg/G$ is an isomorphism except
for the group $Sp_{2n}$ over a base with $2$-torsion. In this case this
morphism is only dominant and we compute it explicitly. We compute the
adjoint quotient in some other classical cases, yielding examples where the
formation of the quotient $\fg\to \fg/G$ commutes, or does not commute, with
base change on $S$.

\section{Introduction}

\indent

Let $G$ be a split semisimple Chevalley group scheme over a
base scheme $S$ and let~$\fg$ be its Lie algebra. The quotient of
$\fg$ by the adjoint action of $G$ in the category of schemes affine over~$S$,
that is to say, the spectrum of the sheaf of $G$-invariant functions
of $\fg$, is traditionally called the {\em adjoint quotient} of $\fg$
and denoted $\fg/G$. Let $T\subset G$ be a maximal torus and $\ft$ its
Lie algebra. There is an induced action of the Weyl group $W=W_T$
on~$\ft$ and the inclusion $\ft\subset\fg$ induces a natural morphism
$\pi:\ft/W\to \fg/G$. In this paper, we call it the {\em Chevalley
  morphism}.

The situation where the base is the spectrum of an algebraically
closed field whose characteristic does not divide the order of the Weyl
group is well documented. In this case $\pi$ is an isomorphism,
as proven by Springer and Steinberg~\cite{SS}. It is known also that
the adjoint quotient is an affine space
(see Chevalley \cite{Ch}, Veldkamp \cite{Ve}, Demazure
\cite{De}). There are counter-examples to these statements when the
characteristic divides the order of the Weyl group.
Another difficulty comes from the fact that we are considering the quotient
$\fg/\Ad(G)$ of the Lie algebra, and not $G/\Int(G)$, and at some
point this derivation
causes some trouble (Steinberg \cite[p.51]{St} was also lead to
the same conclusion).

In this paper, we turn our attention to the integral structure of the adjoint
quotient and the Chevalley morphism, including the characteristics
that divide the order of $W$. In other words we are interested in an
arbitrary base scheme~$S$, and in the behaviour of the previous objects under
base change $S'\to S$. It is not hard to extend the results from simple to
semisimple groups, so for simplicity we restrict to simple Chevalley groups.

Our main result (theorems~\ref{theo_chevalley_dominant} and
\ref{theo_chevalley_iso}) is that in most cases the Chevalley morphism is an
isomorphism, therefore reducing the calculation of $\fg/G$ to the calculation
of a quotient by a finite group:

\bigskip

\noindent {\bf Theorem 1}
{\em Let $G$ be a split simple Chevalley group scheme over a base scheme
  $S$. Then the Chevalley morphism $\pi:\ft/W\to\fg/G$ is schematically
  dominant, and is an isomorphism if $G$ is not isomorphic to $Sp_{2n}$, $n\ge
  1$.}

\bigskip

Note that even when the base is a field, this improves the known results. Our
proof follows a classical strategy. The main new imputs are: over a base
field, a close analysis of the root systems and determination of the
conditions of nonvanishing of the differentials of the roots
(lemma~\ref{racine_non_primitive}), and over a general base, a careful
control of the poles along the singular locus for the relative meromorphic
functions involved in the proof. We treat separately the exceptional case
(theorem~\ref{theo_sp2n}):

\bigskip

\noindent {\bf Theorem 2}
{\em If $G=Sp_{2n}$ then the Chevalley morphism is an isomorphism \iff the
  base has no $2$-torsion. Moreover, over an open affine subscheme
  $\Spec(A)\subset S$, the ring of functions of $\fg/G$ is
$$
A[c_2,c_4,\dots,c_{2n}]
$$
where the functions $c_{2i}$ are the coefficients of the characteristic
polynomial $c_{2i}$. The formation of the adjoint quotient commutes with
arbitrary base change.}

\bigskip

We see that for $G=Sp_{2n}$, the formation of the adjoint quotient commutes
with base change. If this was true for all split simple Chevalley groups, then
we could deduce the main theorem~1 above from the case $S=\Spec(\bZ)$ which is
significantly easier (see corollary~\ref{coro_factorial_ring}). Unfortunately
it is not always so, and in order to see this, we study in
detail the orthogonal groups in types $B$ and $D$. Our main result is (see after
the theorem and subsection \ref{invariant_weyl} for the missing notations):

\bigskip

\noindent {\bf Theorem 3}
{\em If $G=SO_{2n}$ or $G=SO_{2n+1}$ then over an open affine subscheme
$\Spec(A)\subset S$, the ring of functions of $\fg/G$ is the following:

\medskip

{\rm (i)} if $G=SO_{2n}$: $A[c_2,c_4,\dots,c_{2n-2},\pf\,;\,x(\pi_1)^{\epsilon_1}
\dots(\pi_{n-1})^{\epsilon_{n-1}}]$, where $x$ runs through a set of
generators of the $2$-torsion ideal $A[2]\subset A$, and
$\epsilon_i=0$ or $1$, not all $0$,

\medskip

{\rm (ii)} if $G=SO_{2n+1}$:
$A[c_2,c_4,\dots,c_{2n},\,;\,x(\pi_1)^{\epsilon_1}
\dots(\pi_{n})^{\epsilon_{n}}]$, where $x$ runs through a set of
generators of $A[2]$ and $\epsilon_i=0$ or $1$, not all $0$.}

\bigskip

The functions that appear in the preceding theorem are the
coefficients of the characteristic polynomial $c_{2i}$, the Pfaffian
$\pf$ and some functions $\pi_i$ which we call the coefficients of the
Pfaffian polynomial. The functions $c_{2i}$ and $\pf$ are invariant,
but the functions $\pi_i$ are invariant only after multiplication by a
$2$-torsion element. The definition of these objets needs some care, since it
is not always the straightforward definition one would think of.

Using the theorem above, we prove that the formation of the adjoint quotient
for the orthogonal groups commutes with a base change $f:S'\to S$ if and
only if $f^*S[2]=S'[2]$, where $S[2]$ is the closed subscheme defined by the
ideal of $2$-torsion. This holds in particular if $2$ is invertible in
$\cO_S$, or if $2=0$ in $\cO_S$, or if $S'\to S$ is flat. We prove also that
if $S$ is noetherian and connected then the quotient is of finite type
over~$S$, and is flat over $S$ if and only if $S[2]=S$ or $S[2]=\emptyset$.

\bigskip

We feel it useful to say that when we first decided to study the adjoint
quotient over a base other than a field, we started with some examples among
the classical Chevalley groups and considered their Lie algebras.
To our surprise, already in the classical case we could not find concrete
descriptions of them in the existing literature
(for example the Lie algebra of $PSL_n$ over $\Z$).
This lead to our study of the classical Lie algebras over arbitrary bases
(subsection \ref{subsection_classical}). We also faced the problem of
relating the Lie algebra of a group scheme and of any finite quotient of it
(subsection \ref{subsection_quotient}); note that this subsection holds for
any smooth group scheme, not necessarily affine over the base.
Let us finally mention that spin groups over $\bZ$ 
have also been studied very recently in such a concrete
way by Ikai \cite{Ik1}, \cite{Ik2}.

\bigskip

Here is the outline of the article. In the end of this section 1 we give our
notations and prove a combinatorial
lemma about root systems which is crucial in all the paper. In section
\ref{section_Lie_algebras_of_Chev_gps} we give two dual exact sequences
$$
0\to \sLie(K)^\vee\to \sLie(G)^\vee\to \omega^1_{H/S}\to 0
$$
and
$$
0\to \sLie(G)\to \sLie(K)\to (\omega^1_{H/S})^\dag\to 0
$$
describing the relation between the Lie algebra of a smooth group scheme $G$
and the Lie algebra of a quotient $K:=G/H$ (see more precise assumptions in
propositions \ref{exact_sequence_for_Lie} and
\ref{prop_quotient_groupe_fini}). Then we specialize to Chevalley groups and
their Lie algebras over $\bZ$. We describe their weight decomposition
(subsection \ref{subsection_algebre_lie}), the
intermediate quotients of $G\to G^{ad}$ and $\Lie(G)\to \Lie(G^{ad})$
(subsection \ref{diff_quotient_maps}) and we illustrate our results by
describing the classical Chevalley Lie algebras
(\ref{subsection_classical}). In section \ref{section_morphisme_de_Chev} we
prove theorem~1 above. In the remaining
sections \ref{case_S0_2n}, \ref{case_S0_2n+1} and \ref{case_Sp_2n} we treat
the examples of theorems~2 and~3 above by computing explicitly the map
$\ft/W\to \fg/G$ (see theorem~\ref{theo_so2n}, corollary~\ref{coro_so2n},
theorem~\ref{theo_so2n+1}, theorem~\ref{theo_sp2n}).

\tableofcontents

\subsection{General notations} \label{notations}

\indent

All rings are commutative with unit. If $A$ is a ring, we denote by $A[2]$ its
$2$-torsion ideal, defined by $A[2]=\{a\in A,\,2a=0\}$. If $S$ is a scheme, we
denote by $S[2]$ its closed subscheme defined by the $2$-torsion ideal sheaf.

If $X$ is an affine scheme over $\Spec(A)$ we always denote by $A[X]$ its
  function ring.

If $S$ is a scheme, $X$ is a scheme over $S$, and $T\to S$ is a base change
morphism, we denote by $X\times_S T$ or simply $X_T$ the $T$-scheme obtained
by base change. In all the article, we call relative Cartier divisor of $X$
over $S$ an effective Cartier divisor in $X$ which is flat over $S$.

Finally, the linear dual of an $\cO_S$-module $\cF$ is denoted
$\cF^\vee:=\Hom_{\cO_S}(\cF,\cO_S)$.

\subsection{Notations on group schemes}

\indent

Let $S$ be a scheme and let $G$ be a group scheme over $S$.
We will use the following standard notation: $e_G:S\to
G$ is the unit section of $G/S$, $\Omega^1_{G/S}$ is the sheaf of relative
differential $1$-forms of $G/S$, and $\omega^1_{G/S}=e_G^*\Omega^1_{G/S}$.
Recall that $\Omega^1_{G/S}=f^*\omega^1_{G/S}$ where $f:G\to S$ is the
structure map, so that $\Omega^1_{G/S}$ is locally free over $G$ if and only
if $\omega^1_{G/S}$ is locally free over $S$.

We will write $\Lie(G/S)$ (or simply $\Lie(G)$) for the Lie algebra of $G/S$,
and $\sLie(G/S)$ (or simply $\sLie(G)$) for the sheaf of sections of
$\Lie(G/S)$. Note that $\Lie(G/S)$ is the vector bundle $\bV(\omega^1_{G/S})$,
with the Grothendieck notation. Sometimes we shall also use gothic style
letters for Lie algebras, like $\fg$, $\ft$, $\fpsl$, $\fso$, etc.

By Chevalley group scheme over a scheme $S$, we mean a deployable reductive
group scheme over $S$, with the terminology of
\cite[expos{\'e} XXII, d{\'e}finition 1.13]{sga}.
By \cite[expos{\'e} XXIII, corollaire 5.3]{sga},
such a group is characterized up to isomorphism by its type
(as defined in \cite[expos{\'e} XXII, d{\'e}finition 2.7]{sga}: this is 
essentially the root datum together with a module included in the weight 
lattice and containing the root lattice)
and is equal to $G_S$, where $G$ is a Chevalley
group scheme over the ring of integers.


\subsection{Roots that are integer multiples of weights}

The next lemma comes up at various places in the article. It has as
a consequence the fact that the differential
of a root can vanish along the Lie
algebra of a maximal torus in a simple Chevalley group only in case
this group is $Sp_{2n}$ (including $Sp_2 \simeq SL_2$) - see
lemma~\ref{lemm_racine}. This will be crucial throughout the article:
lemma~\ref{lemm_birational}, on which relies the proof of
theorem~\ref{theo_chevalley_iso}, is again a consequence of this lemma, as well
as the fact that a simply-connected Lie algebra is equal to its own derived
algebra in all cases but $\fsp_{2n}$, see proposition~\ref{prop_commutateur}.

\begin{lemm}    \label{racine_non_primitive}
Let $R$ be a simple reduced root system, $Q(R)$ the root lattice and $P(R)$
the weight lattice. Assume there exists 
$\alpha \in R,\lambda \in P(R),l \in \N$ such that 
$\alpha = l.\lambda$ and $l \geq 2$.
Then $l=2$, and either $R$ is of type $A_1$,
or $R$ is of type $C_n$ and $\alpha$
is a long root.
\end{lemm}
\begin{proo}
Let us assume that $R$ is the root system defined in
\cite[Planches I to IX]{bourbaki}.
If $R$ is of type $A_1$, then the roots are
$\alpha = \epsilon_1-\epsilon_2$ and $-\alpha$. Since
$\epsilon_1 - (\epsilon_1+\epsilon_2)/2 = (\epsilon_1-\epsilon_2)/2$ 
is a weight, $\alpha$ is indeed twice a weight. Now let's assume that
$R$ is of rank greater than 1.

The hypothesis of the lemma implies that
\begin{equation}     \label{scalaire_l}
\forall \beta \in R, \scal{\beta^\vee , \alpha} 
= l\scal{\beta^\vee , \lambda} \in l \Z.
\end{equation}
Let $\beta$ be a root. If $\alpha$ and $\beta$ have the same lentgh,
by \cite[VI, no 1.3, Proposition 8]{bourbaki},
$\scal{\beta^\vee,\alpha} \in \{-1,0,1\}$.
If moreover we know that $\scal{\beta^\vee,\alpha} \not = 0$, we see that
(\ref{scalaire_l}) cannot hold.
By \cite[VI, no 1, proposition 15, p. 154]{bourbaki}, we can assume
that $\alpha$ is a simple root.
This implies that in the Dynkin diagram of
$R$, all edges containing the vertex corresponding to $\alpha$ must be
multiple edges.

This excludes all simply-laced root systems, as
well as the root system of type $F_4$.
Moreover, if $R$ is of type $B_n$ with $n \geq 3$, 
then $\alpha$ has to equal $\alpha_n$, but since
$\scal{\alpha_{n-1}^\vee,\alpha_n} = -1$, we have a contradiction.
If $R$ is of type $C_n$, then $\alpha$
has to equal $\alpha_n$ again. Since $\alpha_n = 2\epsilon_n$, we have
indeed $\alpha \in l.P(R)$ with $l=2$.
Since $B_2=C_2$, the last case to be settled is that of $G_2$. But in
this case $Q(R) = P(R)$; since $R$ is reduced, it is not
possible that a root be a multiple of a weight.
\end{proo}




\section{On the Lie algebra of Chevalley groups}
\label{section_Lie_algebras_of_Chev_gps}

\subsection{Lie algebras of quotients and coverings}
\label{subsection_quotient}

In this section, our aim is to relate the Lie algebra of a group $G$
and the Lie algebra of a quotient $G/H$.
More precisely, we consider a scheme $S$, a flat $S$-group
scheme $G$, and a subgroup scheme $H\subset G$ which is flat and of finite
presentation over $S$. It follows from a theorem of Mickael Artin
\cite[corollary 6.3]{artin} that the
quotient fppf sheaf $K:=G/H$ is representable by an algebraic space over
$S$. In the cases of interest to us, it will always be representable by a
scheme. We let $\pi\colon G\to K$ denote the quotient morphism, and
$e_K:=\pi\circ e_G$. Whether or not $H$ is normal, we write $\Lie(K)$ for the
restriction of the tangent space along $e_K$ and $\sLie(K)$ for its sheaf of
sections.

\begin{prop}
\label{exact_sequence_for_Lie}
Assume that $G,H,K$ are as above.
\begin{trivlist}
\itemn{1} There is a canonical exact sequence of quasi-coherent
$\cO_S$-modules:
$$
\omega^1_{K/S}\to \omega^1_{G/S}\to \omega^1_{H/S}\to 0\;,
$$
where $\omega^1_{G/S}\to \omega^1_{H/S}$ is the
natural map deduced from the inclusion $H \subset G$.
\itemn{2} Assume furthermore that $G$ is smooth over $S$ and that
there is schematically dominant morphism $i:U\to S$ such that
$H\times_S U$ is smooth over $U$. Then, there is a canonical exact
sequence of coherent $\cO_S$-modules:
$$
0\to \sLie(K)^\vee\to \sLie(G)^\vee\to \omega^1_{H/S}\to 0
$$
and $\sLie(G)^\vee\to \omega^1_{H/S}$ is the natural map
deduced from the inclusion $H \subset G$.
\end{trivlist}
\end{prop}

Typically, in the applications, $U$ will be an open subscheme of $S$
or the spectrum of the local ring of a generic point.

\begin{proo}
(1) We have the fundamental exact sequence for differential $1$-forms:
$$
\pi^*\Omega^1_{K/S}\to \Omega^1_{G/S}\to \Omega^1_{G/K}\to 0\;.
$$
By right-exactness of the tensor product,
the sequence remains exact after we pullback via $e_G$. The only thing left to
prove is that there is a canonical isomorphism $e_G^*\Omega^1_{G/K}\simeq
\omega^1_{H/S}$. In order to do so, we use the fact that $G\to K$ is an
$H$-torsor, so that we have an isomorphism $t\colon H\times_S G \to
G\times_K G$ given by $t(h,g)=(hg,g)$. We consider the fiber square:
$$
\xymatrix{H\times_S G \ar[r]^t & G\times_K G \ar[r]^{\quad\pr_2} \ar[d]_{\pr_1}
  & G \ar[d]^\pi \\
& G \ar[r]^\pi & K}
$$
Then, if we call $f\colon H\times_S G\to H$ the projection, we have the
sequence of isomorphisms on $H\times_S G$:
$$
t^*\pr_1^*\Omega^1_{G/K}\;\simeq\; t^*\Omega^1_{G\times_K G/G}\;\simeq\;
\Omega^1_{H\times_S G/G}\;\simeq\; f^*\Omega^1_{H/S}
$$
(the first and the third isomorphisms come from the invariance of the module
of relative diffe\-rentials by base change, \cite{EGA},~IV.16.4.5).
Pulling back along $e_H\times e_G$, we get the desired result. Moreover,
following the identifications, we see that the map
$\omega^1_{G/S}\to \omega^1_{H/S}$ is the same as the map induced
by the inclusion $H \subset G$.

\medskip

\noindent (2)
Since $G$ is smooth over $S$ and $G\to K$ is faithfully flat, then $K$ is also
smooth over $S$. Hence $\cM=\omega^1_{K/S}$ and $\cN=\omega^1_{G/S}$
are locally free $\cO_S$-modules of finite rank, so that
$$
\cM\simeq \sLie(K)^\vee
\quad \mbox{and} \quad
\cN\simeq \sLie(G)^\vee \;.
$$
It remains to check that $\cM\to\cN$ is injective. This will follow
from the diagram
$$
\xymatrix{i_*i^*\cM\ \ar@{^(->}[r] & i_*i^*\cN \\
\overset{}{\cM} \ar@{^(->}[u] \ar[r] & \overset{}{\cN} \ar[u]\\}
$$
if we describe the injective morphisms therein.
Since $i:U\to S$ is schematically dominant and $\cM$ is flat, we have
an injective morphism $\cM\to \cM\otimes i_*\cO_U$ and the target
module is isomorphic to $i_*i^*\cM$ by the projection formula.
Besides, the morphism $G\times_S U\to {G/H}\times_S U$ is smooth since
$H\times_S U$ is smooth over $U$, so by the short exact sequence of
$\Omega^1$'s for a smooth morphism, the morphism $i^*\cM\to i^*\cN$ is
injective.
By left exactness the morphism $i_*i^*\cM\to i_*i^*\cN$ is
injective also.
\end{proo}

If $H$ is finite over $S$, like in the cases we have in mind, we can dualize
the exact sequence of proposition~\ref{exact_sequence_for_Lie} thanks to
a Pontryagin duality for certain torsion modules, which we now
present. Let~$A$ be a commutative ring and let $Q$ be the total quotient ring
of $A$, i.e. the localization with respect to the multiplicative set of
nonzerodivisors (in fact we should better consider the module of global
sections of the sheaf of total quotient rings on $\Spec(A)$, but in this
informal discussion it does not matter). We wish to associate to any finitely
presented torsion $A$-module $M$ a dual $M^\dag=\Hom_A(M,Q/A)$. For general
$M$ this does not
lead to nice properties such as biduality; for example if $A=k[x,y]$ is a
polynomial ring in two variables and $M=A/(x,y)$ it is easy to see that
$M^\dag=0$. In this example there is a presentation $A^2\to A\to M\to 0$ but
one can see that there is no presentation $A^n\to A^m\to M\to 0$ with
$n=m$. In fact, this is a consequence of our results below. Note that the fact
that $M$ is torsion implies $n\ge m$, thus if we can find a presentation with
$n=m$ it is natural to say that $M$ has {\em few relations}. By the structure
theorem for modules over a principal ideal domain, all finite abelian groups
have few relations, and from our point of view, this is the crucial property
of finite abelian groups that makes Pontryagin duality work. These
considerations explain the following definition.

\begin{defi} \label{defi_TMFR}
Let $\cF$ be a coherent $\cO_S$-module; denote by $\cK$ the sheaf of total
quotient rings of $\cO_S$. We say that $\cF$ is a
{\em torsion module with few relations} if $\cF\otimes\cK=0$ and $\cF$ is
locally the cokernel of a morphism $(\cO_S)^n\to (\cO_S)^n$ for some $n\ge 1$.
\end{defi}

We have the following easy characterization:

\begin{prop} \label{characterization_of_TMFRs}
Let $\varphi:\cE_1\to\cE_2$ be a morphism between locally free $\cO_S$-modules
of the same finite rank and let $\cF=\coker(\varphi)$. Then $\cF$ is a torsion
module with few relations if and only if the sequence
$0\to \cE_1\to\cE_2\to \cF\to 0$ is exact, i.e. $\varphi$ is injective.
\end{prop}

\begin{proo}
If $\varphi$ is injective, then locally over an open set where $\cE_1$ and
$\cE_2$ are free, its determinant $\det(\varphi)\in\cO_S$ is a
nonzerodivisor. Therefore
$\varphi\otimes\Id:\cE_1\otimes\cK\to\cE_2\otimes\cK$ is surjective, hence an
isomorphism. It follows that $\cF\otimes\cK=\coker(\varphi\otimes\Id)=0$.
Conversely if $\cF$ is a torsion module with few relations, then
$\coker(\varphi\otimes\Id)=\cF\otimes\cK=0$ so that $\varphi\otimes\Id$ is an
isomorphism. Since $\cE_1$ and $\cE_2$ are flat we have injections
$$
\xymatrix{\cE_1\otimes\cK\ \ar@{^(->}[r] & \cE_2\otimes\cK \\
\overset{}{\cE_1} \ar@{^(->}[u] \ar[r] & \overset{}{\cE_2} \ar@{^(->}[u] \\}
$$
and it follows that $\varphi$ is injective.
\end{proo}

\begin{defi} \label{PD}
Given a coherent $\cO_S$-module $\cF$ we define its {\em Pontryagin dual} by
$$\cF^\dag=\cH om_{\cO_S}(\cF,\cK/\cO_S)\;.$$
\end{defi}

As the following proposition proves,
there is a satisfactory duality if we restrict to torsion modules with few
relations.

\begin{prop} \label{Pontryagin_duality}
Let $\cF$ be a torsion $\cO_S$-module with few relations. Then:
\begin{trivlist}
\itemn{1} $\cF^\dag$ is also a torsion $\cO_S$-module with few relations, and
the canonical morphism $\cF\to \cF^{\dag\dag}$ is an isomorphism.

\itemn{2} For each exact sequence $0\to \cE_1\to\cE_2\to \cF\to 0$ where
$\cE_1$, $\cE_2$ are locally free $\cO_S$-modules of the same finite rank, we
have a canonical dual exact sequence $0\to \cE_2^\vee\to\cE_1^\vee\to \cF^\dag\to
0$.
\end{trivlist}
\end{prop}

\begin{proo}
The assertions in point (1) are local over $S$ and therefore are easy consequences
of point~(2). In order to prove point (2) we set
$\cG=\coker(\cE_2^\vee\to\cE_1^\vee)$ and we construct a canonical
nondegenerate pairing $\cF\times\cG\to \cK/\cO_S$, as follows. Since $\cF$ is
torsion and finitely generated, locally (over an
open subset $U\subset S$) there is a nonzerodivisor $a\in\cO_S$ such that
$a\cE_2\subset\cE_1$. Given two sections $f\in\cE_2$ and $g\in\cE_1^\vee$ over
$U$, we let $\langle f,g\rangle$ denote the class of $\frac{1}{a}g(af)\in \cK$
modulo $\cO_S$. It is easy to check that this is independent of the choice
of~$a$. If $f\in\cE_1$ or if $g\in\cE_2^\vee$ then $\langle f,g\rangle=0$ so
there results a pairing $\cF\times \cG\to \cK/\cO_S$ and we will now check
that it induces isomorphisms $\cF\to\cG^\dag$ and $\cG\to\cF^\dag$. By
symmetry we will consider only $\sigma:\cG\to\cF^\dag$. If $\langle
\cdot,g\rangle$ is zero then we claim that $g\in \cE_1^\vee$ extends to a form
on $\cE_2$. Indeed, for each $f\in\cE_2$ we have $\frac{1}{a}g(af)\in \cO_S$
so that the definition $g(f):=\frac{1}{a}g(af)$ is unambiguous, since $a$ is a
nonzerodivisor. It follows that $\sigma$ is injective. In order to check
surjectivity we may assume that $S$ is the spectrum of a local ring, and in
this case $\cE_1$, $\cE_2$ are trivial. Any $u:\cF\to \cK/\cO_S$ factors
through $\frac 1a\cO_S/\cO_S\subset \cK/\cO_S$ and then induces a morphism
$\cE_2\to \cO_S/a\cO_S$. Since $\cE_2$ is trivial this map lifts to
$u':\cE_2\to\cO_S$. Moreover if $x \in \cE_1$ then $u'(x) \in a \cO_S$,
so we can set $v(x) = \frac 1a u'(x)$; 
then it is easy to check that
$v$ is a form $g$ on $\cE_1$ that gives rise to $u$. Hence
$\sigma$ is surjective.
\end{proo}

If $S$ is a Dedekind scheme, that is to say a n\oe therian normal scheme of
dimension~$1$, then all coherent torsion $\cO_S$-modules are torsion modules
with few relations (by the structure theorem for modules of finite
type). However in general it is not so, as soon as $\dim(S)\ge 2$, and we saw
a counter-example before definition \ref{defi_TMFR}.

We are now able to dualize the sequence of Lie
algebras~\ref{exact_sequence_for_Lie}~(2) either if $H$ is smooth or
if it is finite.

\begin{prop} 
\label{prop_quotient_groupe_fini}
Let $G$ be a smooth $S$-group scheme and $H\subset G$ a subgroup scheme
which is flat and of finite presentation over $S$. Let $K=G/H$ be the
quotient.
\begin{trivlist}
\itemn{1} If $H$ is smooth over $S$, then we have an exact
sequence of locally free Lie algebra $\cO_S$-modules
$$0\to\sLie(H)\to\sLie(G)\to\sLie(K)\to 0$$
and if furthermore $K$ is commutative, we have
$$[\sLie(G),\sLie(G)]\subset \sLie(H)\;.$$
\itemn{2} If $H$ is finite over $S$ and there is a
schematically dominant morphism $i:U\to S$ such that $H\times_S U$ is
{\'e}tale over $U$, then there is a canonical exact sequence of coherent
$\cO_S$-modules:
$$
0\to \sLie(G)\to \sLie(K)\to (\omega^1_{H/S})^\dag\to 0 \;,
$$
and if furthermore $H$ is commutative, then
$$[\sLie(K),\sLie(K)]\subset \sLie(G)\;.$$
\end{trivlist}
\end{prop}

\begin{proo}
(1) All the sheaves in the exact
sequence~\ref{exact_sequence_for_Lie}~(2) are locally free, so
dualization yields the asserted result. It is clear that the resulting
sequence is an exact sequence of sheaves of Lie algebras, so
$[\sLie(G),\sLie(G)]\subset \sLie(H)$ in case $K$ is commutative.

\noindent (2) The exact sequence~\ref{exact_sequence_for_Lie}~(2)
and proposition~\ref{characterization_of_TMFRs} imply that
$\omega^1_{H/S}$ is a torsion module with few relations. We get the dual
sequence from proposition~\ref{Pontryagin_duality}. Here
$(\omega^1_{H/S})^\dag$ is not a Lie algebra, so it is a little more subtle to
deduce that $[\sLie(K),\sLie(K)]\subset \sLie(G)$.
The assertion is local on $S$ so we may assume that
$H$ is embedded into an abelian scheme $A/S$
(that is to say a smooth proper group scheme over $S$ with geometrically
connected fibers), by a theorem of Raynaud
(\cite{BBM}, Theorem 3.1.1). Let $\pi:A\to B=A/H$ be the quotient abelian
scheme, and let $G'=(G\times_S A)/H$ where $H$ acts by
$h(g,a)=(hg,h^{-1}a)$. We have two exact sequences of smooth
$S$-schemes:
$$
1\to G\to G'\stackrel{p}{\to} B\to 1
$$
and
$$
1\to A\to G'\to K\to 1 \;.
$$
By smoothness we derive exact sequences of sheaves of Lie algebras
$$
0\to \sLie(G)\to \sLie(G')\stackrel{p}{\to} \sLie(B)\to 0
$$
and
$$
0\to \sLie(A)\stackrel{i}{\to} \sLie(G')\to \sLie(K)\to 0\;.
$$
Combining these exact sequences we have an exact sequence
$$
0\to \sLie(G)\to\sLie(K)\to\sLie(B)/\pi(\sLie(A))\to 0
$$
where $\pi=p\circ i$. Here, the arrow
$\sLie(K)\to\sLie(B)/\pi(\sLie(A))$ is induced by $p$ which is a
morphism of Lie algebras. It follows immediately that
$[\sLie(K),\sLie(K)]\subset \sLie(G)$.
\end{proo}

\subsection{Lie algebras of Chevalley group schemes}

\label{subsection_algebre_lie}

\indent

Let $G$ be a split simple Chevalley group scheme over $\Z$, $T \subset G$ a
split maximal torus over $\bZ$, and write as in \cite{sga} $T=D_\bZ(M)$, where
$M$ is a free $\bZ$-module. Since $G$ is smooth, $\Lie(G)$ is a vector bundle
and hence is determined by $\sLie(G)$. Since the base is affine, this is in
turn determined by the free $\Z$-module $\sLie(G)(\Z)=\Lie(G)(\Z)$ together
with its Lie bracket.

\begin{prop}        \label{prop_decomposition}
There is a weight decomposition
$$\Lie(G)(\Z) = \Lie(T)(\Z) \oplus \bigoplus_\alpha \Lie(G)(\Z)_\alpha$$
over the integers. Moreover, letting $Q(R)$ (resp. $P(R)$)
denote the root (resp. weight) lattice, we have
$Q(R) \subset M\subset P(R)$.
\end{prop}

\begin{proo}
Since $G$ is a smooth split reductive group scheme over $\Z$, this essentially
follows from \cite{sga}, Expos{\'e}~I, 4.7.3, as explained in \cite{sga},
Expos{\'e}~XIX, no 3.
\end{proo}

Now let $H$ be a subgroup of the center of $G$ and let $i_H$ denote the
inclusion of the character group of $T/H$ in that of $T$.

\begin{prop}
Under the natural inclusions 
$$\Lie(G)(\Z) \subset \Lie(G)(\Q) = \Lie(G/H)(\Q) \supset \Lie(G/H)(\Z) \ ,$$
we have $\Lie(G)(\Z)_{i_H(\alpha)} = \Lie(G/H)(\Z)_\alpha$.
\end{prop}
\begin{proo}
Since $H\subset T$, by proposition \ref{prop_quotient_groupe_fini}, there are
injections $\Lie(G)(\Z) \subset \Lie(G/H)(\Z)$ and $\Lie(T)(\Z) \subset
\Lie(T/H)(\Z)$, both of index $|H|$. All these maps are compatible with the
injection in $\Lie(G)(\Q)$. Thus for each $\alpha$, the inclusion
$\Lie(G)(\Z)_{i_H(\alpha)} \subset \Lie(G/H)(\Z)_\alpha$ must be of index 1,
proving the proposition.
\end{proo}

\begin{rema}   \label{rema_serre}
The Lie algebra over $\Z$ defined by generators and relations by Serre
\cite{serre} is the simply-connected one, that is to say the Lie algebra of
the simply-connected corresponding group scheme, because, with his
notations, the generators $H_i$ are by definition the coroots.
\end{rema}

Recall that $\pi : G \rightarrow G/H$ denotes the quotient morphism.

\begin{prop}
\label{prop_commutateur}
Assume $G$ is simply-connected.
\begin{trivlist}
\itemn{1} When $G$ is not $Sp_{2n}$, $n\ge 1$, we have
$$[\Lie(G)(\Z),\Lie(G)(\Z)] = \Lie(G)(\Z)$$
and $[\Lie(G/H)(\Z),\Lie(G/H)(\Z)] = d\pi(\Lie(G)(\Z))$.
\itemn{2}
When $G = Sp_{2n}$, then $[\Lie(G)(\Z),\Lie(G)(\Z)]$ has index $2^{2n}$ in
$\Lie(G)(\Z)$.
\end{trivlist}
\end{prop}

\begin{proo}
(1) Let $\g = \Lie(G)(\Z)$ and
choose a Cartan $\Z$-subalgebra
$\fh \subset \g$. Choose a basis of the roots, and denote by
$\fu_+,\fu_- \subset \g$
the direct sum of the positive (resp. negative) root-spaces.
By corollary \ref{prop_decomposition}, we have
$\g = \h \oplus \u_+ \oplus \u_-$.
Since $G$ is neither $SL_2 (=Sp_2)$ nor $Sp_{2n}$, by 
lemma~\ref{racine_non_primitive},
no root is an integer multiple of a weight, and so 
$[\fh,\fu_\pm] = \fu_\pm$. Moreover, it follows from Serre's presentation of
the simple Lie algebras in terms of the Cartan matrix
(see remark~\ref{rema_serre}) that in this case
$[\g,\g] \supset \h$.
In particular
$$
\begin{array}{rcl}
d\pi(\Lie(G)(\Z)) & = & d\pi([\Lie(G)(\Z),\Lie(G)(\Z)]) \\
& = & [d\pi(\Lie(G)(\Z)),d\pi(\Lie(G)(\Z))] \\
& \subset & [\Lie(G/H)(\Z) , \Lie(G/H)(\Z)] \ . \\
\end{array}
$$
The reverse inclusion follows from proposition
\ref{prop_quotient_groupe_fini}.

\lpara

\noindent (2)
Assume that $G$ stabilises the form $\matdd 0{I_n}{-I_n}0$, where $I_n$
stands for the identity matrix. Then
$$
\g = \left \{ \matdd AB C{-\tr A} : \tr B = B,\tr C = C \right \}.
$$
If $A$ is an arbitrary matrix and $B$ is symmetric, then we have the
equality
$$\left [ \matdd A0 0{-\tr A}, \matddr 0B 00 \right ] 
= \matdd 0{AB+B\tr A} 00.$$
From this it follows that
$[\g,\g] \subset \g$ is the $\Z$-submodule of elements
$\matdd  AB C{-\tr A}$ with $B$ and $C$ having even diagonal elements.
Therefore, it is a submodule of index $2^{2n}$.
\end{proo}

\subsection{The differential of the quotient maps}
\label{diff_quotient_maps}

We will now describe the differentials of the quotient maps between Chevalley
groups in the neighbourhood of a prime $p\in\Spec(\bZ)$. So we consider the
base ring $R=\bZ_{(p)}$. Let $G$ be simply-connected and let $n$ be the order
of the center of $G$. Assume moreover that the center of $G$ is the group of
$n$-th roots of unity $\mu_{n}$ (this is the case if $G$ is not of type
$D_{2l}$; for this particular case see subsection
\ref{subsection_classical}). Write $n=p^km$ with $m$ prime to $p$, and
$G_i:=G/\mu_{p^i}$. We have the successive quotients
$$
G=G_0\to G_1\to G_2\to\dots\to G_k\to G^{ad}
$$
and the corresponding sequence of Lie algebras
$$
\Lie(G)=\Lie(G_0)\to \Lie(G_1)\to\Lie(G_2)\to\dots\to\Lie(G_k)
\stackrel{\simeq}{\longrightarrow} \Lie(G^{ad}) \ .
$$
On the generic fibre all these maps are isomorphisms. In order to study what
happens on the closed fibre, we set $\fg=\Lie(G)(\F_p)$ and
$\fg_i=\Lie(G_i)(\F_p)$, and we let $\fz_i$ resp. $\fz$ denote the center of
$\fg_i$ resp. $\fg$. We start with a lemma:

\begin{lemm}       \label{lemm_centre}
The center $\fz_i$ is isomorphic to the one-dimensional Lie algebra $\F_p$
if $i<k$, and the algebra $\fg_k$ has trivial center.
\end{lemm}

\begin{proo}
Let $x\in\fg_i$ be a central element.
According to the decomposition of proposition \ref{prop_decomposition}, we can
write $ x = \sum x_\alpha + h$. 
The lemma is easily checked directly
when $\fg_i=\fsl_2$ or $\fg_i=\fsp_{2n}$, so
assume we are not in these cases.

According to the following lemma~\ref{lemm_racine}, for
any root $\beta$ we may assume
that there exists $t \in \ft$ such that $d\beta (t) \not = 0$.
We then have $0 = [ t , x ] = \sum d\alpha(t) x_\alpha$, from which it
follows that $x_\beta = 0$. Thus $x=h \in \ft$.
Now, let again $\beta$ be an arbitrary root and let $0 \not = y \in
(\fg_k)_\beta$. We have
$0 = [ x , y ] = d\beta(x) . y$, therefore $d\beta(x) = 0$.

Since we can reverse the above argument, the center of $\fg_i$ consists of
all the elements in $\ft$ along which all the roots vanish.
With the notations of proposition \ref{prop_decomposition},
$\ft \simeq M^\vee \otimes \F_p$, and therefore
$\fz_i \simeq \Hom(M/Q(R),\F_p)$. Since $M/Q(R) \subset P(R)/Q(R)$ and
in our case $P(R)/Q(R)$ is principal, $M/Q(R)$ is also principal and
$\fz_i$ can be at most 1-dimensional. Moreover, it is trivial \iff $Q(R)=M$,
which means that $\fg_i$ is adjoint, or $i=k$.
\end{proo}

\begin{lemm}      \label{lemm_racine}
Assume that $\fg$ is neither isomorphic to $\fsl_2$ nor $\fsp_{2n}$, or
that the characteristic of $k$ is not $2$. Then there exists a finite extension
$K$ of $k$ and $t \in \ft \otimes K$ such that
$\forall \alpha \in R$, $d\alpha(t) \not = 0$.
\end{lemm}

\begin{proo}
Let $R$ denote the root system of $G$ and let $p$ denote
the characteristic of $k$; we have by proposition
\ref{prop_decomposition} $Q(R) \subset M \subset P(R)$.
The linear functions on $\ft$ defined over $k$ are in bijection with $M
\otimes k$; therefore a root $\alpha \in Q(R)$ will yield a vanishing function
on $\ft$ \iff it is a $p$-multiple of some element in $M$. By lemma
\ref{racine_non_primitive}, this can occur only if $p=2$, $M = P(R)$ (thus $G$
is simply-connected) and $G$ is of type $A_1$ or $C_r$. By assumption we are
not in these cases. Taking a finite extension $K/k$ if needed, $\fg$ is not a
union of a finite number of hyperplanes, so the lemma is proved.
\end{proo}

Let $\fg' := \fg / \fz$.
We can now describe the maps $\Lie(G_i)\to\Lie(G_{i+1})$ on the closed fibre:

\begin{prop}
The Lie algebras $\fg_i$ are described as follows:
\begin{trivlist}
\itemn{1}
for all $i$ with $0 < i < k$, we have an isomorphism of Lie algebras
$\fg_i\simeq \fg' \oplus \F_p$.
In particular, all these Lie algebras are isomorphic.
\itemn{2}
for $i=0$ we have a non-split exact sequence of Lie algebras
$0 \rightarrow \F_p \rightarrow \fg_0 
\rightarrow \fg' \rightarrow 0$.
\itemn{3}
for $i=k$ we have a non-split exact sequence of Lie algebras
$0 \rightarrow \fg' \rightarrow \fg_k
\rightarrow \F_p \rightarrow 0$.
\end{trivlist}
In these terms, the maps $\fg_i\to \fg_{i+1}$ are described as follows.
The map $\fg_0\to \fg_1$ takes $\F_p$ to zero and maps onto
$\fg'\subset \fg_1$, and for all $i$ with $0 < i < k$, the map $\fg_i\to
\fg_{i+1}$ takes $\F_p$ to zero and maps $\fg'\subset \fg_i$ isomorphically
onto $\fg'\subset \fg_{i+1}$.
\end{prop}

\begin{proo}
Let $Z_i:=\ker(G_i\to G_{i+1})$. For all $i\le k-1$, we have
$Z_i\simeq \mu_p$, and its Lie algebra is included in $\fz_i$. Because
$Z_i\to G_{i+1}$ is trivial, the map $\fg_i\to \fg_{i+1}$ takes
$\fz_i$ to $0$.

By tensoring the result of proposition \ref{prop_quotient_groupe_fini}
by $\F_p$, there is an exact sequence
$$
\fg_i \to \fg_{i+1} \to \bZ/p\bZ \to 0\ ,
$$
from which it follows that $\fg_i/\fz_i$ is mapped isomorphically onto a
codimension $1$ subalgebra of $\fg_{i+1,\bF_p}$ denoted
$\fg'_{i+1,\bF_p}$. By lemma \ref{lemm_centre}, no
$x \in \fg'_{i+1,\bF_p}$ can be central in $\fg_{i+1}$ so that we have, for
$0 < i < k$,
$$\fg_{i,\F_p} = \fg_{i,\F_p}' \oplus \fz_{i,\F_p}$$
as vector spaces. Since $\g_{i,\F_p}'$ is a Lie subalgebra, it is also an
equality of Lie algebras. In particular all the Lie algebras
$\fg_i$ for $0 < i < k$ are isomorphic.

For $i=0$, we have an exact sequence of Lie algebras
$0 \rightarrow \F_p \rightarrow \fg_{0,\F_p} 
\rightarrow \fg' \rightarrow 0$, but this
sequence does not split (in fact if it did split, then we would have
$[\fg_{0,\F_p} , \fg_{0,\F_p}] \subset \fg'$, contradicting proposition
\ref{prop_commutateur}).

For $i=k>0$, we have a sequence
$0 \rightarrow \fg' \rightarrow \fg_{k,\F_p} 
\rightarrow \F_p \rightarrow 0$, which again
cannot split because $\fg_{k,\F_p}$ has trivial center,
by lemma \ref{lemm_centre}.
\end{proo}

\begin{rema}
For $0<i<k$, consider the Lie algebras $\Lie(G_i)$, as schemes over
$\Spec(\bZ_{(p)})$. They have isomorphic underlying vector bundles (namely the
trivial vector bundle), and their generic fibres as well as their special
fibres are isomorphic as Lie algebras. However, they need not be
isomorphic. For example, if $G=SL_{p^k}$ then it is immediate from
proposition~\ref{prop_commutateur} that $G_i=G/\mu_{p^i}$ uniquely determines
$i$, because $[\Lie(G_i)(\bZ),\Lie(G_i)(\bZ)]=\Lie(G)(\bZ)$ and the quotient
$\Lie(G_i)(\bZ)/\Lie(G)(\bZ)$ is a cyclic abelian group of order
$p^i$. In fact the only difference between them comes from the definition of
the Lie bracket on the total space.
\end{rema}

\subsection{Classical Lie algebras}

\label{subsection_classical}

We now give an explicit description of some of the classical Chevalley Lie
algebras, in which the above sequence of Lie algebras will become very
transparent.

Let $M$ be a free $\Z$-module of rank $n$, and let $m$ be an integer
dividing $n$. We define a $\Z$-Lie algebra $L(M|m)$ as follows:
let $\Hom(M,\frac 1m M)$ denote the $\Z$-module
of linear maps $f:M\otimes\bQ\to M\otimes\bQ$ such that $f(M)\subset \frac 1m
M$. Any such map induces a map $\overline f:M/mM\to \frac 1m M/M$.
Note that multiplication by $m$ induces a canonical isomorphism $\frac 1m
M/M\simeq M/mM$ so that $\overline f$ may be seen as an endomorphism of the
free $\bZ/m\bZ$-module $M/mM$.
Finally, let $L(M|m)$ (resp. $S(M|m)$) denote the submodule of $\Hom(M,\frac 1m
M)$ of elements $f$ such that $\overline f$ is a homothety
(resp. a homothety with vanishing trace). These are obviously
Lie subalgebras of $\End(M\otimes\bQ)$.

\begin{prop}
\label{prop_psln}
Let $n,m$ be as above; then $\Lie(SL_n / \mu_m)(\Z) \simeq S(\Z^n|m)$.
\end{prop}

\begin{proo}
Let $n,m$ be integers such that $m$ divides $n$. Let $\fsl_n$ denote
the Lie algebra of $SL_n$, and let $\fsl_{n,m}$ denote the Lie algebra
of the quotient $SL_n/\mu_m$.
Obviously, we have
$\fsl_n(\Z) \subset \fsl_n(\Q)$,
$\fsl_{n,m}(\Z) \subset \fsl_{n,m}(\Q)$, and
$\fsl_n(\Q) = \fsl_{n,m}(\Q)$ is the usual Lie algebra over $\Q$ of
traceless matrices.

The exact sequence of proposition \ref{exact_sequence_for_Lie}
translates in our case to
$$
0 \rightarrow \fsl_{n,m}(\Z)^\vee \rightarrow \fsl_n(\Z)^\vee \rightarrow
\Z/m\Z \rightarrow 0.
$$
Note that $\fsl_n(\Z)^\vee \simeq \fgl_n(\Z)^\vee / (trace)$, so that
evaluation at $I$ modulo $m$ is well-defined and yields the last arrow.
Let $f_{i,j} : \fsl_n(\Q) \rightarrow \Q$ be the linear form
$M \mapsto M_{i,j}$.
Thus $\fsl_{n,m}(\Z)^\vee \subset \fsl_n(\Z)^\vee$ is the set of linear forms
$\sum \lambda_{i,j}f_{i,j}$ with $m|\sum \lambda_{i,i}$. Since
$\fsl_{n,m}(\Z)$ is the dual in $\fsl_n(\Q)$ of $\fsl_{n,m}(\Z)^\vee$, it
follows applying $f_{i,j}$ $(i \not = j)$ (resp.
$m.f_{i,i}$, $f_{i,i} - f_{j,j}$) that if $M \in \fsl_n(\Q)$ belongs to
$\fsl_{n,m}(\Z)$, then $M_{i,j} \in \Z$ (resp. 
$m.M_{i,i} \in \Z, M_{i,i}-M_{j,j} \in \Z$). Conversely, matrices satisfying
these conditions are certainly in $\fsl_{n,m}(\Z)$.
The proposition is therefore proved.
\end{proo}

Now, assuming that $n$ is even, we describe the four Lie algebras of type $D_n$:
$\fspin_{2n}$, $\fso_{2n}$, $\fpso_{2n}$, and a fourth one that we denote by
$\fpspin_{2n}$. We need to introduce some notation.

\begin{nota}
Let $n$ be an integer and let us consider the following sublattices of $\Q^n$:
\begin{itemize}
\item
$N_{z} = \bZ^n$. 
\item
$N_{ad}$ is generated by $\Z^n$ and $\frac 12 (1,\ldots,1)$.
\item
If $l$ is even, $N_{ps}$ is the sublattice of $N_{ad}$ of elements $(x_i)$
with $\sum x_i$ divisible by 2.
\item
$N_{sc}$ is the sublattice of $\Z^n$ of elements $(x_i)$ with $\sum x_i$
  divisible by 2.
\item
We denote by $L_*$ ($* \in \{z,ad,ps,sc \}$)
the lattice of matrices with off-diagonal coefficients in $\bZ$
and with the diagonal in $N_*$.
\end{itemize}
\end{nota}

\begin{prop}  
\begin{trivlist}  
\itemn{1} There is a natural identification of  
$\fsp_{2n}(\bZ)$ (resp. $\fpsp_{2n}(\bZ)$) with the Lie algebra  
of matrices of the form  
$\matdd ABC{-\tr A}$ with $\tr B=B,\tr C=C$ $(n\times n)$-matrices
with coefficients in $\bZ$ and  
$A \in \fgl_n(\bZ)$ (resp. $ A \in L(\bZ^n|2)$).  
\itemn{2} Assume $n$ is odd.
Then there is a natural identification of $\fso_{2n}(\bZ)$  
(resp. $\fpso_{2n}(\bZ)$, $\fspin_{2n}(\bZ)$) with the Lie algebra of  
matrices of the form $\matdd ABC{-\tr A}$ with $\tr B = -B, \tr C = -C$  
and $A$ in $L_z$ (resp. $L_{ad}$, $L_{sc}$).
\itemn{3} Assume $n$ is even. There is a natural identification of
$\fso_{2n}(\bZ)$  
(resp. $\fpso_{2n}(\bZ)$, $\fspin_{2n}(\bZ)$, $\fpspin_{2n}(\bZ)$)
with the Lie algebra of  
matrices of the form $\matdd ABC{-\tr A}$ with $\tr B = -B, \tr C = -C$  
and $A$ in $L_z$ (resp. $L_{ad}$, $L_{sc}$, $L_{ps}$).
\end{trivlist}  
\end{prop}  
  
\begin{proo}
This is a direct consequence of proposition \ref{prop_decomposition}. For
example, let $\fg(\bZ)$ be a Lie algebra of type $D_n$ over $\bZ$.
Proposition \ref{prop_decomposition} implies 
that all Lie algebras of type $D_n$ will differ only
by their Cartan subalgebras. Therefore $\fg(\bZ)$ is in fact a set of matrices
of the form $\matdd ABC{-\tr A}$ with $\tr B = -B, \tr C = -C$, and the 
off-diagonal coefficients of $A$ in $\bZ$. Moreover, with the notations of
proposition \ref{prop_decomposition}, $\fg(\bZ)$ corresponds to a module $M$ between
the root lattice and the weight lattice, and the Cartan subalgebra (the subalgebra
when $A$ is diagonal and $B=C=0$)
identifies with the dual lattice of $M$. Therefore the 
description of the proposition follows from the description of the root and
weight lattices in \cite{bourbaki}.
\end{proo}

Using this description or lemma \ref{lemm_centre} we can deduce
the dimension of the center of $\fg(\F_2)$:

\medskip

$$
\begin{array}{lcccc}
& \fspin_{2n}(\F_2) & \fso_{2n}(\F_2) & \fpspin_{2n}(\F_2) & \fpso_{2n}(\F_2) \\
n \mbox{ even } & 2 & 1 & 1 & 0 \\
n \mbox{ odd } & 1 & 1 & - & 0 \\
\end{array}
$$

\medskip

For example, we give a description of the center of $\fspin_{2n}(\F_2)$.
Let $C_1$ (resp. $C_2$) be the matrix of the form $\matdd ABC{-\tr A}$
with $B=C=0$ and $A=I_n$ (the identity matrix, resp. $A=\matddr 2000$, the
matrix with only one non-vanishing coefficient in the top-left corner, equal
to 2). Note that $C_1,C_2 \in \fspin_{2n}(\bZ)$ (but $C_1,C_2$ are not
divisible by 2 in $\fspin_{2n}(\bZ)$). For any matrix $B$ in $\fspin(\bZ)
\subset \fso(\bZ)$, we obviously have $[C_1,B]=0$ and $2 | [C_2,B]$. From
these remarks it follows that the classes of $C_1$ and $C_2$ modulo
$2.\fspin(\bZ)$ generate over $\F_2$ the 2-dimensional center of
$\fspin_{2n}(\F_2)$.




\section{The Chevalley morphism}

\label{section_morphisme_de_Chev}

We start this section by some elementary results on regular elements in Lie
algebras of algebraic groups (\ref{section_regular_elements}), to be used
in~\ref{subsection_Chevalley_iso}, then we prove that the Chevalley morphism
is unconditionnally schematically dominant
(\ref{subsection_chevalley_dominant}) and finally we prove the the Chevalley morphism
is an isomorphism unless $G=Sp_{2n}$ for some $n\ge 1$
(\ref{subsection_Chevalley_iso}). The case $G=Sp_{2n}$ will be treated later.

\subsection{Regular elements} \label{section_regular_elements}

Let $\fg$ be a restricted Lie algebra of dimension $d$ over a field. For each
$x\in \fg$, we denote by $\chi(x)=t^d+c_1(x)t^{d-1}+\dots+c_{d-1}(x)t+c_d(x)$
the characteristic polynomial of $\ad x$ acting on~$\fg$. The {\em rank}
of $\fg$ is the least integer $l$ such that $c_{d-l}\ne 0$, and we set
$\delta:=c_{d-l}$. An element $x\in \fg$ is called {\em regular} if the
nilspace $\fg_0(\ad x):=\ker (\ad x)^d$ of $\fg$ relative to $\ad x$ has
minimal dimension~$l$.
An element $x$ is regular if and only if $\delta(x)\ne 0$.
Note that our definition of regular elements differs from that in \cite{Ve},
according to which the singular locus has codimension 3 in $\fg$.
The Cartan
subalgebras of minimal dimension are exactly the centralizers of regular
elements. For these facts see \cite{Str}, pages 52-53.

If, furthermore, $\fg$ is the Lie algebra of a smooth connected group $G$,
then in fact the Cartan subalgebras are conjugate, and in particular they all
have the same dimension. This is proven in \cite{Fa}, corollary 4.4. Also, in
this case the coefficients of the characteristic polynomial $\chi$ are
invariant functions for the adjoint action of $G$ on $\fg$.

We will use the notation $\Sing(\fg)$ for the closed subscheme of singular
elements, defined by the equation $\delta=0$, and $\Reg(\fg)$ for its
complement, the open subscheme of regular elements.
We have the corresponding subschemes $\Sing(\ft)$ and $\Reg(\ft)$ in $\ft$.

In the relative situation, if $\fg$ is a Lie algebra of dimension $d$ over a
scheme $S$, then the objects $\chi$, $\delta$, $\Reg(\fg)$, $\Sing(\fg)$
are defined by the same procedure as above. We recall our general convention
that a {\em relative Cartier divisor} of some $S$-scheme $X$ is an effective
Cartier divisor in $X$ which is flat over $S$.

\begin{lemm}
\label{lemm_sing_cartier}
Let $G$ be a split simple Chevalley group over a scheme $S$, not isomorphic to
$Sp_{2n}$, $n\ge 1$. Let $s:\Sing(\fg)\to S$ be the locus of singular
elements. Then $s_*\cO_{\Sing(\fg)}$ is a free $\cO_S$-module, in particular
$\Sing(\fg)$ is a relative Cartier divisor of $\fg$ over $S$.
\end{lemm}

\begin{proo}
Since the objects involved have formation compatible with base change, it is
enough to prove the lemma over $S=\Spec(\bZ)$.
We have to prove that the ring $\bZ[\fg]/(\delta)$ is free as a
$\bZ$-module. Since $\delta$ is homogeneous, this ring is graded.
If we can prove that it is flat over $\bZ$, then its homogeneous components
are flat also, and since they are finitely generated, they are free over
$\bZ$, and the result follows. So it is enough to prove that
$\bZ[\fg]/(\delta)$ is flat. By the corollary to theorem 22.6 in \cite{Ma}, it
is enough to prove that the coefficients of $\delta$ generate the unit ideal,
or in other words, that $\delta$ is a nonzero function modulo each prime
$p$. So we may now assume that the base is a field $k$ of characteristic $p\ge
0$, and we may also assume that $k$ is algebraically closed.
Let $\ft$ be the Lie algebra of a maximal torus $T$.
By lemma~\ref{lemm_racine}, we can choose
$t \in \ft_k$ such that $\forall \alpha \in R$, $d\alpha(t) \not = 0$.
Then $\delta(t)$ is the product of the scalars $d\alpha(t)$, up to a
sign. Hence it is nonzero.
\end{proo}

We continue with the split simple Chevalley group $G$ over $S$. In the sequel,
products are understood to be fibred products over $S$.
We now turn our attention to the morphism $G/T\times\ft\to \fg$.
We use the same construction as in \cite[3.17]{SS}: note that the normalizer
$N_G(T)$ acts on $G\times \ft$ by $n.(g,\tau) = (gn^{-1},\Ad(n).\tau)$ and
this induces an action of $W=N_G(T)/T$ on $G/T \times \ft$. The morphism
$G/T\times \ft\to \fg$ induced by the adjoint action is clearly $W$-invariant.

\begin{lemm} \label{lemm_birational}
Let $G$ be a split simple Chevalley group over a scheme $S$, not isomorphic to
$Sp_{2n}$, $n\ge 1$. Then the map $G/T\times \ft\to \fg$ is schematically
dominant. Its restriction
$$
b:G/T\times \Reg(\ft)\to \Reg(\fg)
$$
is a $W$-torsor and hence induces an isomorphism $(G/T\times \Reg(\ft))/W\to
\Reg(\fg)$.
\end{lemm}

\begin{proo}
Here again, the objects involved have formation compatible with base change,
so it is enough to prove the lemma over $S=\Spec(\bZ)$. We will first prove
that $b$ is a $W$-torsor. It is enough to prove that $b$ is surjective,
\'etale, and that $W$ is simply transitive in the fibres. Indeed, if $b$ is
\'etale then the action must also be free and $b$ induces an isomorphism
$(G/T\times \Reg(\ft))/W\to \Reg(\fg)$.

The map $c=\ad:G\times \Reg(\ft)\to \Reg(\fg)$ is surjective because if $x\in
\Reg(\fg)$, then its centralizer $\fz(x)$ is a Cartan subalgebra, and since
Cartan subalgebras are conjugate, there exists $g\in G$ such that $(\ad
g)(\ft)=\fz(x)$. Thus there is $y\in\ft$ such that $(\ad g)(y)=x$ and clearly
$y$ is regular. We now prove that $c$ is smooth. Since its source and its
target are smooth over $S$, it is enough to prove that for all $s\in S$, the
map $c_s$ is smooth. By homogeneity, it is enough to prove that the
differential of $c_s$ at any point $(1,t)$ with $t \in \Reg(\ft)$ is
surjective.

Then $T_1G_k = \fg_k$ and the tangent map
$\psi=dc : T_1G_k \times \ft_k \rightarrow \fg_k$ is given by 
$(x,\tau) \mapsto [x,t] + \tau$.
Recall that $\fg_k = \bigoplus_{\alpha \in R} \fg_\alpha \oplus \ft_k$,
where $[\tau,x] = d\alpha(\tau)x$ for all $x \in \fg_\alpha$.
Again by lemma~\ref{lemm_racine}, we can choose
$\tau \in \ft_k$ such that $\forall \alpha \in R$, $d\alpha(\tau) \not = 0$. Thus,
we have $\psi(\fg_k \times \{0\}) = \bigoplus_{\alpha \in R} \fg_\alpha$.
Since $\psi(\{0\} \times \ft_k) = \ft_k$, $\psi$ is a surjective linear map.
It follows that $b$ is also surjective and smooth, hence \'etale, by dimension
reasons.

Finally, let $(g,x)$ and $(h,y)$ have the same image in $\fg$, for
$x,y\in\ft$. This means that $(\ad w)(x)=y$, where $w=h^{-1}g$. Thus
$(\ad w)(\fz(x))=\fz((\ad w)(x)=\fz(y)$, that is to say $(\ad
w)(\ft)=\ft$ since $\fz(x)=\fz(y)=\ft$. By \cite[13.2,13.3]{humphreys},
$T$ is the only maximal torus with Lie
algebra $\ft$, so it follows that $w$ normalizes $T$.
Hence $w$ defines an element of the Weyl group $W$.

Now we consider the map $G/T\times \ft\to \fg$. From the preceding discussion it
is dominant in the fibres, and since $G/T\times\ft$ is flat over $S$, the map
is itself schematically dominant by \cite{EGA}, th{\'e}or{\`e}me 11.10.9.
This concludes the proof of the lemma.
\end{proo}

\subsection{The Chevalley morphism is always dominant}

\label{subsection_chevalley_dominant}

We now deal with the cases that are not covered by
lemma~\ref{lemm_sing_cartier}.

\begin{nota}
Consider the following subalgebras of $\fsl_2$
and $\fsp_{2n}$:
\begin{itemize}
\item
Let $\fb \subset \fsl_2$ be the subalgebra of upper-triangular
matrices. 
\item
Let $L$ denote the set of long roots of $\fsp_{2n}$; if $\alpha$ is a
root, denote by $\fsp_{2n,\alpha}$ the corresponding root space.
\item
Let $\fh \subset \fsp_{2n}$ be the sum
$\ft \oplus \bigoplus_{\alpha \in L} \fsp_{2n,\alpha}$.
\end{itemize}
\end{nota}

\begin{lemm}       
\label{lemm_dominant_sur_fibres_sp}
Let $k$ be a field. Then the maps $(SL_2)_k \times \fb_k \rightarrow
(\fsl_2)_k$ and $(Sp_{2n})_k \times \fh_k \rightarrow (\fsp_{2n})_k$, given by
restricting the adjoint action, are dominant. Moreover, $\fh$ is isomorphic,
as a Lie algebra, to $\fsl_2^{\oplus n}$.
\end{lemm}

Combining the two statements of this lemma, it follows that in the proof of
theorem~\ref{theo_chevalley_dominant} below, we will be able to replace the
Lie subaglebra $\fh \simeq \fsl_2^{\oplus n}$ by a sum of the form
$\fb^{\oplus n}$.

\begin{proo}
The result about $(\fsl_2)_k$ is an immediate consequence of the fact
that, over an algebraically closed field, any matrix is conjugated
to an upper-triangular matrix.

To prove that $(Sp_{2n})_k \times \fh_k \rightarrow (\fsp_{2n})_k$
is dominant, we argue as in lemma \ref{lemm_birational}. Since
all the short roots are not integer multiples of a weight, we can
choose $t \in \ft_k$ such that for all short roots $\alpha$,
we have $d\alpha(t) \not = 0$. Let $S$ denote the set of short roots of
$\fsp_{2n}$. If $\psi$ denotes the differential at
$(1,t)$ of the ajoint action, it follows that
$\psi(\fsp_{2n} \times \{0\}) \supset \bigoplus_{\alpha \in S}
\fsp_{2n,\alpha}$. Since $\fh_k = \ft_k \oplus \bigoplus_{\alpha \in L}
\fsp_{2n,\alpha}$, it follows that $\psi$ is surjective and the restriction
of the action is dominant.

To prove that $\h \simeq \fsl_2^{\oplus n}$, one can compute
explicitly in the Lie algebra $\fsp_{2n}$. Assume that $\fsp_{2n}$ is
defined by the matrix $\matdd 0I{-I}0$, where $I$ denotes the identity
matrix of size~$n$. Then a matrix $\matdd ABCD$ belongs to $\fsp_{2n}$
\iff $D = -\tr A$ and $B$ and $C$ are symmetric matrices. Choosing the
torus $\ft = \left \{ \matdd d00{-d} \right \}$ 
in $\fsl_{2n}$, and $\epsilon_i$
the coordinate forms on $\ft$, it is well-known and easy to check that
the long roots are $\pm 2\epsilon_i$. It follows that
$\fh = \left \{ \matdd d \delta \epsilon {-d} : d,\delta,\epsilon
\mbox{ diagonal} \right \}$. Thus $\fh$ is isomorphic to
$\fsl_2^{\oplus n}$.
\end{proo}

\begin{lemm}  \label{lemm_injectivite_sl2}
Let $S=\Spec(A)$ be an affine base scheme and $\g = \fsl_2$.
Then the restriction morphism $A[\fb]^T \rightarrow A[\ft]$ is injective.
\end{lemm}
\begin{proo}
Let $f \in A[\fb]$. Writing a typical element in $\fb$ as
$\matddr ab0a$, we identify $f$ with a polynomial in $a,b$. Since
$$\matddr t00{t^{-1}}\matddr ab0a \matddr t00{t^{-1}} ^{-1} 
= \matdd a{t^2b}0a,$$
$f$ is $T$-invariant \iff $f(a,b) = f(a,t^2b) \in A[a,b,t]$. This means that
$f$ does not depend on $b$, and we indeed have an injection
$A[\fb]^T \rightarrow A[\ft]$.
\end{proo}

Using the two preceding lemmas, we can now prove the main result of this
section:

\begin{theo}     \label{theo_chevalley_dominant}
Let $S$ be a scheme and let $G$ be a split simple Chevalley group over~$S$.
Then the Chevalley morphism $\pi:\ft/W\to\fg/G$ is schematically dominant.
\end{theo}

\begin{proo}
First, let $S=\Spec(\bZ)$. Let us write $\fh = \fb$ in the case of $SL_2$,
$\fh = \fsl_2^{\oplus n}$ in the case of $Sp_{2n}$, and $\fh = \ft$ in the
other cases. We also write $H = T$ (resp. $H=SL_2^{\oplus n} , H=T$) in the
case of $SL_2$ (resp. $Sp_{2n}$, the other cases). The adjoint action
restricts to a map  $\varphi:G \times \fh \to \fg$. If $G$ acts on itself by
left translation, trivially on $\fh$, and by the adjoint action on $\fg$, then
$\varphi$ is $G$-equivariant. Moreover, by lemmas \ref{lemm_birational} and
\ref{lemm_dominant_sur_fibres_sp}, the restriction $\varphi_k$ of $\varphi$ to
any fiber of $\Spec(\Z)$ is dominant. Since the schemes $G \times \fh$ and
$\fg$ are flat over $\Z$, it follows from \cite[th{\'e}or{\`e}me
  11.10.9]{EGA} that $\varphi$ is universally schematically dominant.

Now we let $S$ be arbitrary, and prove that $\pi$ is schematically dominant.
The question is local over $S$ so we may assume $S=\Spec(A)$ affine. 
On the function rings, the Chevalley morphism can be decomposed
as two successive restriction morphisms $A[\fg]^G  \to A[\fh]^H \to A[\ft]^W$.

First we concentrate on $A[\fg]^G  \to A[\fh]^H$. We have already proven that
the map $\varphi^* : A[\fg] \rightarrow A[G] \otimes_A A[\fh]$ is
injective. This map is $G$-equivariant, so if $f \in A[\fg]$ is $G$-invariant,
we have
$$
\varphi^*(f) \in
(A[G] \otimes_A A[\fh])^G = A[G]^G \otimes_A A[\fh]
= A \otimes_A A[\fh] = A[\fh] \ .
$$
Therefore $\varphi^*(f) = 1 \otimes i^*(f)$ where $i:\fh \rightarrow \fg$ is
the inclusion. Since $\varphi^* : A[\fg] \rightarrow A[G] \otimes_A A[\fh]$ is
injective, it follows that $A[\fg]^G \rightarrow A[\fh]$ is injective.

In case $G \not = Sp_{2l}$, we have $\ft = \fh$ so the proof of the
theorem is complete. In case $G = SL_2$, lemma \ref{lemm_injectivite_sl2}
shows the injectivity of the Chevalley morphism. Finally, let us
consider the case of $Sp_{2n}$ with $n>1$. Since by lemma
\ref{lemm_dominant_sur_fibres_sp}, $\fh \simeq \fsl_2^{\oplus n}$
and $H \simeq SL_2^n$, and since the theorem is proved for $SL_2$, we know
that $A[\fh]^H  \to  A[\ft]$ is injective, and we can once again conclude.
\end{proo}

It is interesting to mention an easy consequence of this theorem~: if the base
scheme is $\Spec(\bZ)$, then the Chevalley morphism is an isomorphism, see the
corollary below. This will of course be a particular case of
theorems~\ref{theo_chevalley_iso} and \ref{theo_sp2n}, however, while
theorem~\ref{theo_chevalley_iso} needs some more work, we get the present
result almost for free. This corollary is in fact worthwhile because it shows
that if the formation of the adjoint quotient commuted with base change (which
is not the case), then we would get the case of a general base scheme $S$ from
the case $S=\Spec(\bZ)$ and hence everything would be finished right now.
So here is this corollary:

\begin{coro} \label{coro_factorial_ring}
Assume that $S$ is the spectrum of a factorial ring with
characteristic prime to the order of $W$. Then the Chevalley morphism
$\pi:\ft/W\to\fg/G$ is an isomorphism.
\end{coro}

\begin{proo}
Let $S=\Spec(A)$ and let $K$ be the fraction field of $A$. By
theorem~\ref{theo_chevalley_dominant} it is enough to prove that the
restriction morphism $\res:A[\fg]^G \to A[\ft]^W$ is surjective. Let $P$ be
a $W$-invariant function on $\ft$. From the assumption on the
characteristic of $K$, it follows that $K[\fg]^G \to K[\ft]^W$ is an
isomorphism, so there is $Q\in K[\fg]^G$ such that $\res(Q)=P$. Since
$A$ is factorial, we can write $Q=cQ_0$ where $Q_0$ is a primitive
polynomial (i.e. the gcd of its coefficients is $1$) and $c\in K$ is the
content of $Q$. If we write $c=r/s$ with $r$ and $s$ coprime in $A$,
we claim that $s$ is a unit in $A$. For, otherwise, some prime $p\in
A$ divides $s$. Then $\res(\bar{rQ_0})=\bar{sP}=0$ in
$(A/p)[\ft]^W$ so $\res(\bar{Q_0})=0$ in
$(A/p)[\ft]^W$, since $r$ and $s$ are coprime. By
theorem~\ref{theo_chevalley_dominant} again, it follows that $\bar{Q_0}=0$ in
$(A/p)[\fg]^G$, in contradiction with the fact that $Q_0$ is
primitive. Hence $Q\in A[\fg]^G$ as was to be proved.
\end{proo}

\subsection{The Chevalley morphism is an isomorphism for $G\ne Sp_{2n}$}

\label{subsection_Chevalley_iso}

Before proving theorem~\ref{theo_chevalley_iso}, our main result,
we state two lemmas which will ultimately allow us to show that a function
vanishes on the locus of singular elements.

\smallskip

The first lemma is a general result about the ``indicator function'' of the
locus of the base where a function has a zero along a fixed divisor~:

\begin{lemm} \label{indicator_of_zeroes}
Let $X\to S$ be a morphism of schemes and let $D\subset X$ be a relative
effective Cartier divisor with structure morphism $p:D\to S$, such that
$p_*\cO_D$ is a locally free $\cO_S$-module. Let $*$ be the
one-point set, and given a global function $f$ on $X$, let $F$ be the
functor on the category of $S$-schemes defined as follows~: for any $S$-scheme
$T$,
$$
F(T)=\left\{
\begin{array}{cl}
* & \mbox{if $f_T$ has a zero along $D_T$}, \\
\emptyset & \mbox{otherwise.} \\
\end{array}
\right.
$$
Then $F$ is representable by a subscheme $S_0\subset S$ which is open and
closed in $S$.
\end{lemm}

\begin{proo}
We first prove that $F$ is represented by a closed subscheme $S_0\subset S$.
The assertion is local on $S$ so we may assume $S$ affine and $p_*\cO_D$
free. Let $f_i$ be the finitely many components of $f_{|D}$ on some basis of
$\Gamma(D,\cO_D)$ as a free $\Gamma(S,\cO_S)$-module. Then, $f$ has a zero
along $D$ if and only if all
the $f_i$ vanish. Since $D$ is a relative Cartier divisor, the formation of
these objects commutes with base change, so that the above description is
functorial. Then obviously $F$ is represented by the closed subscheme of
finite presentation $S_0\subset S$ defined by the ideal of $\Gamma(S,\cO_S)$
generated by the coefficients $f_i$.

It order to prove that $S_0\subset S$ is also open, we will use another
description of $F$. Let $L=\Spec(\Sym(\cO_X(D))$ and $d:L\to\bA^1_X$ be the
morphism of line bundles over $X$ induced by the canonical morphism of
invertible sheaves $\cO_X\to\cO_X(D)$. We can view $f$ as a section of the
trivial line bundle, i.e. $f:X\to\bA^1_X$. Then to say that $f$ has a zero
along $D$ means that there exists a section $\varphi:X\to L$ of $L$ such that
the diagram
$$
\xymatrix{X \ar[r]^\varphi \ar[rd]_f & L \ar[d]^d \\
& \bA^1_X \\}
$$
commutes. Since $D$ is a relative Cartier divisor, such a $\varphi$ is unique
if it exists. In other words, $F(T)$ is the set of sections $\varphi:X_T\to
L_T$ such that $f_T=d_T\circ\varphi$, or again, if we set $Z:=d^{-1}(f(X))$,
then $F(T)$ is the set of sections $\varphi:X_T\to L_T$ such that $\varphi(X_T)$
is a closed subscheme of $Z_T$. In the sequel, we will take the freedom to
write simply $X$ for the divisor
$\varphi(X)\subset L$, and $D$ for the pullback of $D$ via the structure
morphism $L\to X$. Since $\varphi$ is determined by its image, we have
described $F$ as an open subfunctor of a Hilbert functor of $Z/S$.

To prove that $S_0\subset S$ is open, it is enough to prove that this is a
smooth morphism. Since $S_0\subset S$ is of finite presentation, it is enough
to prove that it is formally smooth at all points $s\in S_0$. To do this, we
will show that the deformations of $\varphi_s$ are unobstructed.
Let $k$ be the residue field of $s$, let $A$ be a local artinian
$\cO_S$-algebra with residue field $k$ and assume that $\varphi_s$ has been
lifted to $\varphi_A$. We have to prove that for each nilthickening $A'\to A$,
i.e. a surjective morphism with kernel $M$ annihilated by the maximal ideal
of~$A'$, the space of obstructions to lifting $\varphi_A$ to $A'$
vanishes. By the theory of the Hilbert functor, the obstruction space is
$\Ext^1_{Z_A}(\cI_A,(\cO_Z/\cI)\otimes_k M)$ where $\cI$ is the ideal sheaf of
$X$ in $Z$; note that $(\cO_{Z_A}/\cI_A)\otimes_k M=(\cO_Z/\cI)\otimes_k
M$. We will compute this group by using an explicit resolution of $\cI$.
In the sequel, we write $X$ for $X_A$, $D$ for $D_A$, etc.

We remark that $Z=X+D$ as a sum of Cartier divisors
of $L$. This is not hard to see and we leave the details to the reader.
Now let $a_D:\cO_L\to \cO_L(D)$ and $a_X:\cO_L\to \cO_L(X)$ be the canonical
morphisms. Consider the sequence
$$
\dots
\stackrel{a_X}{\longto}\cO_L(-2Z)\stackrel{a_D}{\longto}\cO_L(-X-Z)
\stackrel{a_X}{\longto}\cO_L(-Z)\stackrel{a_D}{\longto}\cO_L(-X)
\stackrel{a_X}{\longto}\cO_L \ .
$$
Note that this is not even a complex, but we claim that after restricting to
$Z$, we get a resolution of the ideal sheaf $\cI$~:
$$
\dots
\stackrel{a_X}{\longto}\cO_L(-2Z)|_Z\stackrel{a_D}{\longto}\cO_L(-X-Z)|_Z
\stackrel{a_X}{\longto}\cO_L(-Z)|_Z\stackrel{a_D}{\longto}\cO_L(-X)|_Z
\stackrel{a_X}{\longto}\cI\longto 0 \ .
$$
Indeed the image of $a_X:\cO_L(-X)|_Z\to{\cO_L}|_Z=\cO_Z$ is $\cI$, and
locally $X$ and $D$ have equations $t_X$ and $t_D$ in $\cO_L$, and the
sequence is the sequence of $\cO_L/(t_Xt_D)$-modules which is alternatively
multiplication by $t_X$ and $t_D$. From the fact that $t_X$ and $t_D$ are
nonzerodivisors, the exactness of the sequence follows.

So to compute $\Ext^1_Z(\cI,\cO_Z/\cI\otimes_k M)$ we apply
$\Hom_Z(\cdot,\cO_Z/\cI\otimes_k M)$ and take the first cohomology group of
the resulting complex. In fact $\Hom_Z(\cI,\cO_Z/\cI\otimes_k M)=0$~; this is
the tangent space to the functor $F$. Therefore
$\Ext^1_Z(\cI,\cO_Z/\cI\otimes_k M)$ is equal to
$$
\ker\big(a_D:\Hom_Z(\cO_L(-X)|_Z,\cO_Z/\cI\otimes_k
M)\longto\Hom_Z(\cO_L(-Z)|_Z,\cO_Z/\cI\otimes_k M)\big) \ .
$$
Locally $\cO_L(-X)|_Z\simeq \cO_L(-Z)|_Z\simeq \cO_Z$, and
 $a_D$ takes a map $\sigma:\cO_Z\to\cO_Z/(t_X)\otimes_k M$ to
$t_D\sigma$. Since $D$ is a Cartier divisor in $X$, it follows that $a_D$ is
injective. Consequently
$$\Ext^1_Z(\cI,\cO_Z/\cI\otimes_k M)=0 \ ,
$$
thus the functor $F$ is unobstructed, so $S_0$ is formally smooth at~$s$.
\end{proo}

The second lemma proves a statement which is used in \cite{SS}
(proof of 3.17, point (2)).
However we were not able to understand their proof, due to a
vicious circle in the use of an argument from \cite{St}:

\begin{lemm}     \label{lemm_division}
Let $k$ be a field and assume $S = \Spec(k)$. Let $f,g \in k[\fg]^G$
and assume $f_{|\ft} \divise g_{|\ft}$. Then $f \divise g$.
\end{lemm}

\begin{proo}
We may assume that $k$ is algebraically closed.
First assume that $f$ has no square factors.
Let $x \in \fg$ be such that $f(x)=0$; it is then enough to show that $g(x)=0$.
To this end, since $x$ belongs to a Borel subalgebra of $\fg$
(in fact, the Borel subalgebra are the maximal solvable algebras) and
all the Borel subalgebras are conjugated, we may assume that
$x = \tau + \eta$ with $\tau \in \ft$ and
$\eta =\sum_{\alpha > 0} x_\alpha \in \bigoplus_{\alpha > 0} \fg_\alpha$. Let
$X:{\mathbb G}_m \to T$ be a one-parameter-subgroup corresponding to a coweight
$\omega^\vee$ with $\scal{\omega^\vee,\alpha} > 0$ for all positive
roots $\alpha$. We therefore have
$\Ad(X(t)).x = \tau + \sum_{\alpha > 0} t^{n_\alpha} x_\alpha$, with
$\forall \alpha>0,n_\alpha > 0$. Therefore
the closure of the $G$-orbit through $x$ contains
$\tau$, and $f(x)=f(\tau)=0$. Thus $g(x)=g(\tau)=0$.

Thus the lemma is proved in case $f$ is squarefree. Now let $f$ be
arbitrary. Write $f=f_1f_2$ where $f_1$ is the product of the prime factors of
$f$, with multiplicity $1$. So $f_1$ is squarefree, and since $G$ is connected
and hence has no nontrivial characters, we see that $f_1$ is $G$-invariant.
We have $(f_1)_{|\ft} \divise g_{|\ft}$ so $f_1 \divise g$~:
let us write $g=f_1g_2$. By factoriality of $k[\fg]$ we have
$(f_2)_{|\ft} \divise (g_2)_{|\ft}$, so the lemma follows
by induction on the degree of $f$.
\end{proo} 

We can now prove our main result.

\begin{theo}
\label{theo_chevalley_iso}
Let $S$ be a scheme and let $G$ be a split simple Chevalley group over~$S$.
Assume that $G$ is not isomorphic to $Sp_{2n}$, $n\ge 1$.
Then the Chevalley morphism $\pi:\ft/W\to\fg/G$ is an isomorphism.
\end{theo}

\begin{proo}
By theorem~\ref{theo_chevalley_dominant} it is enough to prove that the map on
functions is surjective. Let $f$ be a $W$-invariant function on $\ft$ and let
$f_1$ be the function on $G/T\times\ft$ defined by $f_1(g,x)=f(x)$. Since it
is $W$-invariant, it induces a function on $(G/T\times\ft)/W$ which we denote
by the letter $f_1$ again. By lemma~\ref{lemm_birational}
the function $h:=f_1\circ
b^{-1}$ is a $G$-invariant relative meromorphic function whose domain of
definition contains $\Reg(\fg)$.
We may write $h=k/\delta^m$ for some function $k$ not
divisible by $\delta$, and some integer $m$. Since $k$ is $G$-invariant on a
schematically dense open subset, it is $G$-invariant. Assume that $m\ge 1$.
Let $s\in S$ be a point.
Since a generic element of $\ft_s$ is regular, $h_s$ is defined as a
rational function on $\ft_s$. By definition of $b$, we moreover have
$(h_s)_{|\ft} = f_s$. It follows that we have
$$
(k_s)_{|\ft} = f_s . (\delta_s)_{|\ft}^m.
$$
Therefore the restriction of $\delta_s$ divides the restriction of
$k_s$. Lemma~\ref{lemm_division} implies that $\delta_s$ divides $k_s$. Since
this is true for all $s\in S$, and $p:\Sing(\fg)\to S$ is a relative Cartier
divisor of $\fg/S$ with $p_*\cO_D$ free over $S$ by lemma~\ref{lemm_sing_cartier}, then
we can apply lemma~\ref{indicator_of_zeroes} to conclude that $\delta$ divides
$k$. This is a contradiction with our assumptions, therefore $h$ is a regular
function extending $f$ to a $G$-invariant function on $\fg$.
\end{proo}

In the remaining sections, we compute explicitly the ring of invariants
in the case where $G$ is one of the groups $SO_{2n},SO_{2n+1}$ or $Sp_{2n}$.




\section{The orthogonal group $SO_{2n}$}

\label{case_S0_2n}

In the case of the group $G=SO_{2n}$, the explicit computation will prove that
the formation of the adjoint quotient for the Lie algebra does not commute
with all base changes. In fact, we will be able to describe exactly when
commutation holds. We will see also that over a base field, the quotient is
always an affine space.

\subsection{Definition of $SO_{2n}$} 

\label{subsection_def_SO_2n}

\begin{noth} {\bf The orthogonal group.}
The free $\bZ$-module of rank $2n$ is denoted by $E$; we think of it
as the trivial vector bundle over $\Spec(\bZ)$. The standard quadratic
form of $E$ is defined for $v=(x_1,y_1,\dots,x_n,y_n)$ by
$$
q(v)=x_1y_1+\dots+x_ny_n \ .
$$
It is nondegenerate in the sense that $\{q=0\}\subset
\bP(E)$ is smooth over $\bZ$. The polarization of $q$ is
$$
\langle v,v'\rangle=q(v+v')-q(v)-q(v')=x_1y'_1+x'_1y_1+\dots+x_ny'_n+x'_ny_n \ .
$$
The orthogonal group $O_{2n}$ is the set of transformations $P\in
GL_{2n}$ that preserve $q$, more precisely, the zero locus of the
morphism $\Psi$ from $GL_{2n}$ to the vector space of quadratic forms
defined by $\Psi(P)=q\circ P-q$. Thus the Lie algebra $\fo_{2n}$
is the subscheme of $\fgl_{2n}$ composed of matrices $M$ such that by
$d\Psi_{\Id}(M)=0$ with
$$
d\Psi_{\Id}(M)(v)=\langle v,Mv\rangle \ .
$$
It is not hard to verify that $\fo_{2n}\subset \fgl_{2n}$ is a
direct summand of the expected dimension, so that $O_{2n}$ is a
smooth group scheme over $\bZ$.

\begin{rema} \label{orth_matrix_preserves_polarization}
Let us denote by $B$ the matrix of the polarization of $q$. Clearly, an
orthogonal matrix $P$ preserves the polarization, and it follows that
${}^tPBP=B$. However, one checks easily that the subgroup $X\subset GL_{2n}$
defined by the equations ${}^tP BP=B$ is not flat over $\bZ$ because its
function ring has $2$-torsion. In fact $O_{2n}$ is the biggest subscheme of
$X$ which is flat over $\bZ$. Accordingly $\Lie(X)\subset\fgl_{2n}$ is
defined by ${}^tMB+BM=0$, and $\fo_{2n}$ is the biggest $\bZ$-flat
subscheme of $\Lie(X)$.
\end{rema}
\end{noth}

\begin{noth} {\bf Dickson's invariant.}
Over any field $k$, it is well-known that $O_{2n}\otimes k$ has two
connected components. In odd characteristic, the determinant takes the
value $1$ on one and $-1$ on the other. In characteristic $2$ the
determinant does not help to separate the connected components.
Instead one usually uses Dickson's invariant $D(P)$ defined, for an
orthogonal matrix $P$, to be $0$ if and only if $P$ acts trivially on
the even part of the center of the Clifford algebra. Equivalently,
$D(P)=0$ if and only if $P$ is a product of an even number of
reflections (there is just one exception; see \cite{Ta},
p. 160). Here is a more modern, base-ring-free way to
consider the determinant and Dickson's invariant altogether:

\begin{lemm}
There is a unique element $\delta\in \bZ[O_{2n}]$ such that
$\det=1+2\delta$.
\end{lemm}

\begin{proo}
Since for any $P\in O_{2n}$, we have $\det(P) \in \{-1,1\}$, the
function $\det-1$ vanishes on the fibre $O_{2n}\otimes \bF_2$. Since
$O_{2n}\otimes \bF_2$ is reduced, $2$ divides $\det-1$, yielding the
existence of $\delta$. It is unique because $O_{2n}$ is flat over
$\bZ$, and in particular has no $2$-torsion.
\end{proo}

Let us introduce the $\bZ$-group scheme
$\cG=\Spec(\bZ[u,\frac{1}{1+2u}])$ with unit $u=0$ and
multiplication $u*v=u+v+2uv$. Its fibre at $2$ is isomorphic to the
additive group while all other fibres are isomorphic to the
multiplicative group. When one passes from $\det$ to $\delta$,
the multiplicativity formula $\det(P_1P_2)=\det(P_1)\det(P_2)$ gives
$\delta(P_1P_2)=\delta(P_1)+\delta(P_2)+2\delta(P_1)\delta(P_2)$.
In other words,

\begin{lemm}
$\delta$ defines a morphism of groups $O_{2n}\to \cG$.
\hfill $\square$
\end{lemm}

The schematic image of $\delta$ is the subgroup of $\cG$ given by $u(u+1)=0$,
isomorphic to the constant $\bZ$-group scheme $\zmod{2}$.

\begin{defi}
We define $SO_{2n}$ as the kernel of $\delta$.
\end{defi}

The group $SO_{2n}$ is smooth over $\bZ$ with connected fibres.
The subgroup $T$ of diagonal matrices in $SO_{2n}$ is a maximal torus, we
denote by $\ft$ its Lie algebra and by $\lambda_i$ its coordinate
functions. Its normalizer $N$ is the subgroup of orthogonal monomial
matrices. The Weyl group $W=N/T$ is the semi-direct product
$(\zmod{2})^{n-1}\rtimes\fS_n$ where $\fS_n$ is the symmetric group on
$n$ letters. It acts on $T$ as follows. The subgroup
$(\zmod{2})^{n-1}$ is generated by the transformations
$\varepsilon_{i,j}$ which take $\lambda_i$ and $\lambda_j$ to their
opposite and leave all other $\lambda_k$ unchanged. The subgroup
$\fS_n$ permutes the $\lambda_i$. The action of $W$ on $\ft$ has
analogous expressions that are immediate to write down.
\end{noth}

\begin{noth} {\bf The Pfaffian.}
Recall that there is a unique function on $\fso_{2n}$, called the {\em
  pfaffian} and denoted $\pf$, such that $\det(M)=(-1)^n(\pf(M))^2$. (The
sign $(-1)^n$ comes from the fact that in our context, the pfaffian is
  $\pf'(BM)$ where $\pf'$ is the usual pfaffian.) Furthermore the pfaffian is
  invariant for the adjoint action of $SO_{2n}$.
\end{noth}

\subsection{Invariants of the Weyl group}

\label{invariant_weyl}

\indent

We denote by $\ft$ the $n$-dimensional affine space with coordinate functions
$X_i$, and by $W$ the group generated by the permutations of the coordinates
and the reflections $\varepsilon_{i,j}$ which map $X_i$ and $X_j$ to
their opposite and leave the other coordinates invariant. We denote by 
$\sigma_k$
the complete elementary symmetric functions in $n$ variables.

\begin{prop}      \label{invariant_weyl_so2n}
Let $A$ be a ring, then $A[\ft]^W$ is generated by $X_1\cdots X_n$,
$\sigma_k(X_i^2)$, and 
$x\sigma_k(X_i)$, where $k<n$ and $x$ runs through the 2-torsion
ideal of $A$.
\end{prop}

\begin{proo}
Let $F$ be a function in $X_1,\dots,X_n$ which is invariant under the Weyl
group. Let us say that a monomial is {\em good} if the exponents of its
variables $X_i$ all have the same parity (this is either a monomial in
the $X_i^2$, or $X_1\cdots X_n$ times a monomial in the $X_i^2$). We say that
it is {\em bad} otherwise. We can write uniquely $F$ as the sum of its good
part and its bad part:
$$
F(X_1,\dots,X_n)
=F_1(X_1^2,\dots,X_n^2,X_1\dots X_n)+F_2(X_1,\dots,X_n) \ .
$$
The group $W$ respects this decomposition, hence $F$ being $W$-invariant,
its good and bad parts also are. In particular they are
$\fS_n$-invariant, so that
$$
F(X_1,\dots,X_n)
=G_1(\sigma_1(X_i^2),\dots,\sigma_{n-1}(X_i^2),X_1\dots X_n)
+G_2(\sigma_1(X_i),\dots,\sigma_n(X_i)) \ .
$$
Letting the $\varepsilon_{i,j}$ act, we see that all coefficients of $G_2$
must be $2$-torsion. The proposition is therefore proved.
\end{proo}

\subsection{Computation of the Chevalley morphism}

\label{subsection_invariants_SO_even}

\indent

In this subsection we will describe explicitly the invariants of
$\fso_{2n}$ under $SO_{2n}$ that correspond to the Weyl group invariants
under theorem~\ref{theo_chevalley_iso},
see theorem~\ref{theo_so2n} below.
The Lie algebra $\fso_{2n}$ has a universal matrix $\sM$ whose
most important attributes are its characteristic polynomial $\chi$ and its
pfaffian $\pf=\pf(\sM)$. In fact $\sM$ and $\chi$ are the restrictions of the
universal matrix of $\fgl_{2n}$ and its characteristic polynomial. From the
equality ${}^t\sM B+B\sM=0$ (see \ref{orth_matrix_preserves_polarization}) it
follows that $\chi$ is an even polynomial, that is to say
$$
\chi(t)=\det(t\Id-\sM)=t^{2n}+c_2t^{2n-2}+\dots+c_{2n}.
$$
The functions $c_{2i}$ are invariants of the adjoint action; note that
$$c_{2n}=\det(\sM)=(-1)^n(\pf(\sM))^2 \ .$$
There are some more invariants coming from characteristic $2$. Indeed, in this
case the polar form is alternating, so homotheties are antisymmetric and we can
define the {\em pfaffian characteristic polynomial} by
$$\pi_{\bF_2}(t)=\pf(t\Id-\sM_{\bF_2}) \ .$$
We have $\chi_{\bF_2}(t)=(\pi_{\bF_2}(t))^2$. Now let us consider one
particular lift of $\pi_{\bF_2}$ to $\bZ$:

\begin{defi} \label{defi_coefs_pfaff_polynomial}
Let $\sigma\colon\bF_2\to \bZ$ be such that $\sigma(0)=0$ and $\sigma(1)=1$.
The polynomial $\pi\in\bZ[\fso_{2n}][t]$ is defined as
$\pi(t)=t^n+\pi_1t^{n-1}+\dots+\pi_{n-1}t+\pi_n$ where $\pi_n:=\pf(\sM)$ and
the other coefficients $\pi_i$ ($1\le i\le n-1$) are the lifts of the
corresponding coefficients of $\pi_{\bF_2}$ via $\sigma$.
\end{defi}

We note that for any ring $A$, the (images of the) elements
$\pi_1,\dots,\pi_{n-1},\pi_n$ are algebraically independent over $A$,
because they restrict on a maximal torus to the functions $\sigma_i(X_j)$, the
symmetric functions in the coordinates, which are themselves algebraically
independent over $A$. We defined the functions $\pi_i$ by arbitrary lifting,
but we can make them somehow universal:

\begin{prop} \label{universal_coefs_of_pfaffian}
Let $\cO$ be the ring $\bZ[X]/(2X)$ and denote by $\tau$ the image of $X$
in~$\cO$. Then $(\cO,\tau)$ is universal among rings with a $2$-torsion
element. Moreover, any monomial function
$\tau(\pi_1)^{\alpha_1}\dots(\pi_{n-1})^{\alpha_{n-1}}$ on
$\fso_{2n,\cO}$ is independent of the choice of the lifts $\pi_i$
and invariant under the adjoint action of $SO_{2n,\cO}$. Finally, for each
$i\in\{1,\dots,n\}$ we have $\tau(\pi_i)^2=\tau c_{2i}$.
\end{prop}

\begin{proo}
The universality statement about $(\cO,\tau)$ means that for any pair $(A,x)$
where $A$ is a ring and $x$ is a $2$-torsion element of $A$, there is a unique
morphism $f:\cO\to A$ such that $f(\tau)=x$. This is obvious. Since $2\tau=0$,
it is clear also that $\tau(\pi_1)^{\alpha_1}\dots(\pi_{n-1})^{\alpha_{n-1}}$
is independent of the choice of $\pi_i$. The fact that this monomial is
invariant comes from the invariance of the pfaffian characteristic polynomial
in characteristic $2$. Finally the equalities $\tau(\pi_i)^2=\tau c_{2i}$ come
from the equalities $(\pi_i)^2=c_{2i}$ in characteristic $2$.
\end{proo}

By proposition \ref{universal_coefs_of_pfaffian},
for any ring $A$ and any $x\in A[2]$, the
quantity $x(\pi_1)^{\alpha_1}\dots(\pi_{n-1})^{\alpha_{n-1}}$ is a
well-defined invariant function on $\fso_{2n,A}$.

\begin{theo}      \label{theo_so2n}
Let $A$ be a ring, $G=SO_{2n,A}$, $\fg = \fso_{2n,A}$. The ring of invariants
$A[\fg]^G$ is
$$
A[c_2,c_4,\dots,c_{2n-2},\pf\,;\,x(\pi_1)^{\epsilon_1}
\dots(\pi_{n-1})^{\epsilon_{n-1}}]
$$
where $x$ runs through a set of generators of the $2$-torsion ideal
$A[2]\subset A$, and $\epsilon_i=0$ or $1$, not all $0$.
\end{theo}

\begin{proo}
By theorem~\ref{theo_chevalley_iso} and proposition
\ref{invariant_weyl_so2n}, we have
$A[\fg]^G=A[\sigma_k(X_i^2)\,;\,X_1\dots X_n\,;\,x\sigma_k(X_i)]$ where as
before the $X_i$ are the coordinate functions on the torus. We now use
the previous proposition. Since $c_{2k}$ restricts on the torus to $\pm
\sigma_k(X_i^2)$, the pfaffian restricts to $X_1\dots X_n$, $x\pi_k$ restricts
to $x\sigma_k(X_i)$, and since $x\pi_i^2 = xc_{2i}$, the theorem is proved.
\end{proo}

The behaviour of the ring of invariants is therefore controlled by the
$2$-torsion. More precisely, for a scheme $S$ let $S[2]$ be the
closed subscheme defined by the ideal of $2$-torsion. If $f:S'\to S$ is a
morphism of schemes, we always have $S'[2]\subset f^*S[2]$. We have:

\begin{coro}    \label{coro_so2n}
\begin{trivlist}  
\itemn{1}
The formation of the quotient in the previous theorem commutes with a base
change $f:S'\to S$ if and only if $f^*S[2]=S'[2]$. This holds in particular if
$2$ is invertible in $\cO_S$, or if $2=0$ in $\cO_S$, or if $S'\to S$ is flat.
\itemn{2}
Assume that $S$ is noetherian and connected. Then the quotient
is of finite type over~$S$, and is flat over $S$ if and only if
$S[2]=S$ or $S[2]=\emptyset$.
\end{trivlist}
\end{coro}

\begin{proo}
First we recall some general facts on the formation of the ring of
invariants for the action of an affine $S$-group scheme $G$ acting on an
affine $S$-scheme $X$. The formation of $X/G=\Spec((\cO_X)^G)$ commutes with
flat base change, and in particular with open immersions. It follows that if
$(S_i)$ is an open covering of $S$ then with obvious notation $X_i/G_i$ is an
open set in $X/G$, and $X/G$ can be obtained by glueing the schemes
$X_i/G_i$. Therefore if $(S'_{i,j})$ is an open covering of $S_i\times_S S'$
for all $i$, the formation of the quotient commutes with the base change
$S'\to S$ if and only if for all $i,j$ the formation of the quotient $X_i/G_i$
commutes with the base change $S'_{i,j}\to S_i$. This reduces the proof to the
case of a base change of affine schemes $S'=\Spec(A')\to S=\Spec(A)$.

Call $B$ (resp. $B'$) the ring of invariants over $A$ (resp. $A'$).
Observe that $B$ inherits a graduation from the graduation of the function
algebra of $\fso_{2n,A}$, and its only homogeneous elements of degree~$1$ are
those of the form $x\pi_1$ with $x\in A[2]$. We proceed to prove (1) and (2).

\noindent (1)
The base change morphism $B\otimes_A A'\to B'$ is
$A'[\underline c,\pf,x\underline\pi^\epsilon]\to A'[\underline
  c,\pf,x'\underline\pi^\epsilon]$ where
$$\underline c=(c_2,c_4,\dots,c_{2n-2}) \ , \
x\underline\pi^\epsilon=x(\pi_1)^{\epsilon_1} \dots(\pi_{n-1})^{\epsilon_{n-1}}
 \mbox{ with } x\in A[2]
$$
and $x'\underline\pi^\epsilon$ is the same quantity with $x'\in A'[2]$. This
map is clearly injective. If it is surjective, then in particular for any
$x'\in A'[2]$ we have $x'\pi_1\in A'[\underline
  c,\pf,x\underline\pi^\epsilon]$. Thus there is $a'\in A'$ and
$x\in A[2]$ such that $x'\pi_1=a'x\pi_1$. Since $\pi_1$ is a nonzerodivisor we
get $x'=a'x$, so $A'[2]$ is the image of $A[2]$. This is exactly the assertion
that $f^*S[2]=S'[2]$. The converse is easy, as well as the particular cases
stated in the lemma.

\noindent (2)
If $A$ is noetherian, $A[2]$ is finitely generated and then $B$ is of finite
type over $A$. Now let $I:=A[2]$. If $I=0$ then $B$ is
a polynomial ring, and this is also the case if $I=A$ because then
$c_2,c_4,\dots,c_{2n-2}$ are polynomials in
$\pf(\sM),\pi_1,\dots,\pi_{n-1}$. It remains to prove that if $B$
is flat over $A$ then $I=0$ or $I=A$. In this case the $2$-torsion ideal of
$B$ is $IB$, as we see from tensoring by~$B$ the exact sequence
$$
0\to A/I\stackrel{\times 2}{\to} A\to A/2\to 0 \ .
$$
So for any $y\in I$, we have $y\pi_1\in B[2]=IB$ hence we may write
$y\pi_1=i_1b_1+\dots+i_rb_r$ with $i_k\in I$ and $b_k\in B$. Let $x_k\pi_1$ be
the degree~$1$ component of $b_k$, then by taking the components of degree~$1$
and using the fact that $\pi_1$ is a nonzerodivisor, we find
$y=i_1x_1+\dots+i_rx_r\in I^2$. Thus $I=I^2$, and if the spectrum of $A$ is
connected, this implies $I=0$ or $I=A$.
\end{proo}

\section{The orthogonal group $SO_{2n+1}$}

\label{case_S0_2n+1}

For $G=SO_{2n+1}$, the computation of the quotient is a little more involved
since using the natural representation of dimension $2n+1$ brings some
trouble, as we explain below. We show which point of view on $SO_{2n+1}$ will
lead to the definition of the correct invariants. Then, the results are
essentially the same as for $G=SO_{2n}$.

\subsection{Definition of $SO_{2n+1}$}

\indent

In this section, $E$ is the free $\bZ$-module $\bZ^{2n+1}$.
Its standard quadratic form $q$ is
$$
q(v)=x_1y_1+\dots+x_ny_n+z^2
$$
where $v=(x_1,y_1,\dots,x_n,y_n,z)$. It is nondegenerate, and its polarization
is
$$
\langle v,v'\rangle=q(v+v')-q(v)-q(v')=x_1y'_1+x'_1y_1+\dots+x_ny'_n+x'_ny_n+2zz' \ .
$$
In contrast with the even dimensional case, in characteristic $2$ the
polarization has a nonzero radical which is the line generated
by the last basis vector of $E\otimes\bF_2$.

Now, let $\tilde E=\bZ^{2n+2}$ with canonical basis
$(e_1,e'_1,\dots,e_{n+1},e'_{n+1})$ and standard quadratic form defined (as in
\ref{subsection_def_SO_2n}) by $q(v)=x_1y_1+\dots+x_{n+1}y_{n+1}$ where
$v=(x_1,y_1,\dots,x_{n+1},y_{n+1})$. We consider the isometric embedding
$i\colon E\hookrightarrow \tilde E$ given by
$$
i(x_1,y_1,\dots,x_n,y_n,z)=(x_1,y_1,\dots,x_n,y_n,z,z) \ .
$$
Since $i$ is an isometry, it is harmless to use the same letter for
$q$ and for $q_{|E}$. The orthogonal subspace of $E$ in $\tilde E$ is
the free rank $1$ submodule generated by the vector
$\varepsilon=e_{n+1}-e'_{n+1}$. Note that the group of transformations of
$(\tilde E,q)$ is the group $O_{2n+2}$ as defined in
\ref{subsection_def_SO_2n}. Then we define $SO_{2n+1}$ as a closed
subgroup of $SO_{2n+2}$ by
$$
SO_{2n+1}=\{\,P\in SO_{2n+2},\, P(\varepsilon)=\varepsilon \,\} \ .
$$
Accordingly, its Lie algebra is
$$
\fso_{2n+1}=\{\,M\in \fso_{2n+2},\, M(\varepsilon)=0 \,\} \ .
$$
It is a simple exercise to verify that $\fso_{2n+1}\subset \fso_{2n+2}$ is a
direct summand of the expected dimension, so that $SO_{2n+1}$ is a smooth
group scheme over $\bZ$.

\begin{rema} \label{another_presentation}
Let $O(q|_E)$ be the group of linear transformations of $E$
that preserve~$q$, and $SO(q|_E)$ the kernel of the Dickson invariant. Since a
matrix $P\in SO_{2n+1}$ preserves the line generated by $\varepsilon$, it
preserves its orthogonal~$E$. This leads to a morphism $SO_{2n+1}\to
SO(q|E)$. However, because of the existence of a one-dimensional radical in
characteristic~$2$, one can see that the fibre $SO(q|E)\otimes\bF_2$ is
nonreduced and its reduced subscheme is the subgroup $H$ of transformations
that act as the identity on the radical. Thus $SO_{2n+1}\to SO(q|E)$ is not an
isomorphism. In fact, one may see that this map realizes $SO_{2n+1}$ as
the dilatation of $SO(q|E)$ with center~$H$. Recall from \cite{BLR}, 3.2
that the dilatation is a map $\pi\colon SO_{2n+1}\to SO(q|E)$ which is
universal for the properties: $SO_{2n+1}$ is $\bZ$-flat and its
special fibre at $2$ is mapped into $H$. It can be checked that the dilatation
is indeed smooth over $\bZ$ and is the Chevalley orthogonal group. In this
formulation, the special orthogonal group is not naturally a group of
matrices. This is why we used another presentation. \hfill $\square$
\end{rema}

The subgroup $T\subset SO_{2n+2}$ of diagonal matrices fixing
$\varepsilon$ is a maximal torus of $SO_{2n+1}$, we denote by $\ft$ its Lie
algebra and by $\lambda_i$ its coordinate functions. Its normalizer $N$
is the subgroup of orthogonal monomial matrices fixing
$\varepsilon$. The Weyl group $W=N/T$ is the semi-direct product
$(\zmod{2})^n\rtimes\fS_n$. It acts on $T$ as follows. The subgroup
$(\zmod{2})^{n-1}$ is generated by the transformations $\varepsilon_i$
which take $\lambda_i$ to its opposite and leave
all other $\lambda_k$ unchanged. The subgroup $\fS_n$ permutes the
$\lambda_i$.

\subsection{Explicit computation of the Chevalley morphism}

\label{subsection_invariants_SO_odd}

\indent

Let $\sM$ be the universal matrix over $\fso_{2n+1}$. Using the
embedding of $\fso_{2n+1}$ into $\fso_{2n+2}$, we define
invariants by restriction from those of $\fso_{2n+2}$ defined in
\ref{subsection_invariants_SO_even}. For example, let us view the universal
matrix $\sM$ as a matrix in $\fso_{2n+2}$. Since $\sM(\varepsilon)=0$, the
determinant of~$\sM$ vanishes and hence its characteristic polynomial in
dimension $2n+2$ is
$$
t^{2n+2}+c_2t^{2n}+\dots+c_{2n}t^2 \ .
$$
We {\em define} the {\em characteristic polynomial of $\sM$} as
$$
\chi(t)=t^{2n+1}+c_2t^{2n-1}+\dots+c_{2n}t \ .
$$
Note that this is not the characteristic polynomial associated to
an actual action on the natural representation of dimension $2n+1$.
Using again the embedding in $\fso_{2n+2}$, we see that in characteristic $2$
we have again a polynomial $\pi(t)$ defined uniquely by the identity
$\chi_{\bF_2}(t)=t(\pi_{\bF_2}(t))^2$. By abuse, we call it again {\em
  pfaffian characteristic polynomial}. We may define lifts of its coefficients
by the same process as in definition \ref{defi_coefs_pfaff_polynomial} and we
obtain a polynomial $\pi(t)=t^n+\pi_1t^{n-1}+\dots+\pi_{n-1}t+\pi_n$ where
$\pi_i\in\bZ[\fso_{2n}]$. As in
subsection~\ref{subsection_invariants_SO_even}, for any ring $A$ the elements
$\pi_1,\dots,\pi_{n-1},\pi_n$ are algebraically independent over $A$. In the
same way as in subsection \ref{subsection_invariants_SO_even}, we prove:

\begin{prop}
Let $(\cO,\tau)$ be the ring defined in proposition
\ref{universal_coefs_of_pfaffian}. Then any monomial function
$\tau(\pi_1)^{\alpha_1}\dots(\pi_n)^{\alpha_n}$ on
$\fso_{2n+1,\cO}$ is independent of the choice of the lifts~$\pi_i$
and invariant under the adjoint action of $SO_{2n+1,\cO}$. Also, for each
$i\in\{1,\dots,n\}$ we have $\tau(\pi_i)^2=\tau c_{2i}$.
\hfill $\square$
\end{prop}

So for any ring $A$ and any $x\in A[2]$, the
quantity $x(\pi_1)^{\alpha_1}\dots(\pi_n)^{\alpha_n}$ is a
well-defined invariant function on $\fso_{2n+1,A}$.
Exactly the same proof as the proof of \ref{theo_so2n} gives:

\begin{theo}      \label{theo_so2n+1}
Let $A$ be a ring and $G=SO_{2n+1,A}$. Then the ring of functions of $\fg/G$ is
$$
A[c_2,c_4,\dots,c_{2n}\,;\,x(\pi_1)^{\epsilon_1}\dots(\pi_n)^{\epsilon_n}]
$$
where $x$ runs through a set of generators of the $2$-torsion ideal
$A[2]\subset A$, and $\epsilon_i=0$ or $1$, not all $0$.
\hfill $\square$
\end{theo}

Finally, all the statements of corollary \ref{coro_so2n} hold also word for
word for $G=SO_{2n+1}$.




\section{The symplectic group $Sp_{2n}$}
\label{case_Sp_2n}

\indent

The computation of the adjoint quotient and of the Chevalley morphism
$\pi:\ft/W\to\fg/G$ for $Sp_{2n}$ requires the preliminary computation of the
corresponding quantity for the group $SL_2$. We also found it interesting
to deal with the case of $PSL_2$.

\subsection{Preliminary cases: $SL_2$ and $PSL_2$}

We denote by $\matddr abc{-a}$ the universal matrix of $\fsl_2$; therefore
$A[a]$ is the ring of functions on $\ft$ over $A$.
The same proof as that of proposition \ref{invariant_weyl_so2n} yields:

\begin{fait}      \label{fait_weyl_sl2}
Let $A$ be a ring, then $A[\ft]^W$ is equal to $A[a^2] \oplus aA[2][a^2]$.
\end{fait}

\noindent
We set $\det(a,b,c) = -a^2 - bc$.
This fact and the next proposition show that we don't have
$\ft/W \simeq \fg/G$:

\begin{prop}            \label{prop_sl2}
Let $A$ be a ring, then $A[\fsl_2]^{SL_2} = A[\det]$.
\end{prop}
\begin{proo}
The action of the diagonal matrix $\matddr u00{u^{-1}}$
on the coordinate functions reads:
\begin{equation}     \label{action_diagonale}
\begin{array}{l}
a \mapsto a\\
b \mapsto u^2b\\
c \mapsto u^{-2}c
\end{array}
\end{equation}

Any invariant polynomial can therefore be written as a polynomial in $a$ and
$bc$.

On the other hand the action of the unipotent element $\matddr 1t01$ reads:
\begin{equation}     \label{action_unipotent}
\begin{array}{l}
a \mapsto a+tc\\
b \mapsto b-2ta-t^2c\\
c \mapsto c
\end{array}
\end{equation}

Assume we have a homogeneous invariant of odd degree $2d+1$.
Since it is a polynomial
in $a$ and $bc$, it can be written as $af(a^2,bc)$, with $f$ homogeneous
of degree $d$. We consider the identity
$$
af(a^2,bc) = (a+tc)f((a+tc)^2,(b-2ta-t^2c)c),
$$
and specialise to $a=0$. We get
$tcf(t^2c^2,bc-t^2c^2)=0$ so $f(t^2c^2,bc-t^2c^2)=0$.
Performing the invertible change of coordinates
$d=b+t^2c$, we therefore get $f(t^2c^2,cd)=0=c^df(t^2c,d)$,
from which it follows that $f=0$.

\medskip

Thus there are no invariants of odd degree and the image of the restriction
morphism is included in $A[a^2]$. Since $\det$ is an invariant, this image is
exactly $A[a^2]$, which implies the proposition.
\end{proo}

We pass to $PSL_2$. By proposition \ref{prop_psln} and its proof,
the coordinate ring of $\fpsl_2$ over $A$ is $A[\alpha , b , c]$, where
$\alpha = 2a$.

\begin{fait}      \label{fait_weyl_psl2}
Let $A$ be a ring, then $A[\ft]^W$ is equal to
$A[\alpha^2] \oplus \alpha A[2][\alpha^2]$.
\end{fait}

\begin{prop}           \label{prop_psl2}
Let $A$ be a ring, then $A[\fpsl_2]^{PSL_2} = A[4\det] + \alpha
. A[2][4\det]$.
\end{prop}

\begin{proo}
We know from theorem \ref{theo_chevalley_dominant}
that $A[\fpsl_2]^{PSL_2}$ injects
into $A[\ft]^W = A[\alpha^2] \oplus \alpha A[2][\alpha^2]$. On the other hand,
$4\det = -\alpha^2 - 4bc$ is certainly an invariant in the coordinate ring,
as well as $x\alpha$, if $x \in A$ is a 2-torsion element, since
by (\ref{action_unipotent}) $\alpha$ is mapped to $\alpha + 2tc$ under the
action of the unipotent element $\matddr 1t01$. Thus the proposition is proved.
\end{proo}

\subsection{Explicit computation of the Chevalley morphism}

\indent

We denote by $\ft$ the $n$-dimensional affine space with coordinate functions
$X_i$, and $W$ the group generated by the permutations of the coordinates and
the reflections $\varepsilon_i$ which map $X_i$ to its opposite and leave the
other coordinates invariant. Recall that $\sigma_k$ denotes
the complete elementary symmetric
functions in $n$ variables. The same proof as for proposition
\ref{invariant_weyl_so2n} yields:

\begin{prop}      \label{invariant_weyl_sp2n}
Let $A$ be a ring, then $A[\ft]^W$ is generated by $\sigma_k(X_i^2)$ and
$x\sigma_k(X_i)$, where $k<n$ and $x$ runs through the 2-torsion ideal of $A$.
\hfill $\square$
\end{prop}

We denote by $E$ the natural representation of $G = Sp_{2n}$, of dimension
$2n$. By definition, we therefore have a morphism $G \rightarrow GL(E)$, which
also induces a morphism $\fg \rightarrow \fgl(E)$.
Let $\sM$ be the universal matrix over $\fgl(E)$, and let $\chi$ be its
characteristic polynomial:
$$
\chi(t)=\det(t\Id-\sM)=t^{2n}-c_1t^{2n-1}+c_2t^{n-2}+\dots+c_{2n} \ .
$$

\begin{theo}      \label{theo_sp2n}
Let $A$ be a ring and $G=Sp_{2n,A}$. Then the morphism $\pi:\ft/W\to\fg/G$
is an isomorphism \iff $A$ has no 2-torsion. Moreover,
the ring of functions of $\fg/G$ is
$$
A[c_2,c_4,\dots,c_{2n}] \ .
$$
The formation of the adjoint quotient $\fg\to\fg/G$
over a scheme $S$ commutes with any base
change $S'\to S$.
\end{theo}

\begin{proo}
Let $G=Sp_{2n,A}$ and $\fg = Lie(G)$. By theorem~\ref{theo_chevalley_dominant}
and
proposition~\ref{invariant_weyl_sp2n}, $A[\fg]^G$ is a subring of
$A[\ft]^W=A[\sigma_k(X_i^2)\,;\,x.\sigma_k(X_i)]$. With the notations of
lemma~\ref{lemm_dominant_sur_fibres_sp}, it is also a subring of the image
of $A[\fh]^H$ in $A[\sigma_k(X_i^2)\,;\,x.\sigma_k(X_i)]$.
The latter is $A[\sigma_k(X_i^2)]$
by proposition~\ref{prop_sl2}. Since $c_{2k}\in A[\fg]^G$ maps to $\pm
\sigma_k(X_i^2)\in A[\ft]^W$, the theorem is proved.
\end{proo}

\bigskip\noindent
Pierre-Emmanuel {\sc Chaput}, Laboratoire de Math{\'e}matiques Jean Leray, 
UMR 6629 du CNRS, UFR Sciences et Techniques,  2 rue de la Houssini{\`e}re, BP
92208, 44322 Nantes cedex 03, France. 

\noindent {\it email}: pierre-emmanuel.chaput@math.univ-nantes.fr

\medskip\noindent
Matthieu {\sc Romagny}, Institut de Math{\'e}matiques, Th{\'e}orie des Nombres,
Universit{\'e} Pierre et Marie Curie, Case 82,
4, place Jussieu,
F-75252 Paris Cedex 05.

\noindent {\it email}: romagny@math.jussieu.fr

\end{document}